\documentclass[11pt]{elsarticle}

\usepackage[utf8]{inputenc}
\usepackage{default}
\usepackage{graphicx}
\usepackage{amsfonts,amsmath,amssymb,amsthm}
\usepackage{color}
\usepackage{bm}
\usepackage{algorithm,algorithmic}

\def\be{\mathbf{e}}

\def\bi{\mathbf{i}}

\def\bk{\mathbf{k}}
\def\bl{\mathbf{l}}

\def\bt{\mathbf{t}}

\def\bv{\mathbf{v}}
\def\bw{\mathbf{w}}
\def\bx{\mathbf{x}}

\def\bz{\mathbf{z}}

\def\bA{\mathbf{A}}

\def\cA{\mathcal{A}}

\def\cI{\mathcal{I}}
\def\cJ{\mathcal{J}}
\def\cK{\mathcal{K}}
\def\cL{\mathcal{L}}

\def\bxi{\boldsymbol{\xi}}

\def\sol{\bz}
\def\solh{\bz_h}
\def\solhn{\bz_{h,n}}

\def\qoi{J(\bxi)}
\def\qoih{J_h(\bxi)}
\def\ptqoih{J_h(\bxi_{\bl,\bi})}
\def\qoihn{J_{h,n}(\bxi)}
\def\ptqoihn{J_{h,n}(\bxi_{\bl,\bi})}
\def\qoinh{J_{h,n}(\bxi)}

\def\dest{\delta(\qoih)}
\def\ptdest{\delta(\ptqoih)}
\def\eest{\varepsilon(\solhn,\bphi_{h,n})}

\def\pteest{\varepsilon(J_{h,n}(\bxi_{\bl,\bi}))}
\def\qoien{J_n^{\varepsilon}}
\def\qoienm{J_{n,m}^{\varepsilon}}

\def\gridpt{\brv_{\bl,\bi}}

\newcommand{\intp}[2]{I_{#1}[#2]}
\newcommand{\dom}{\Gamma_{\bxi}}

\newcommand{\norm}[2]{\left\lVert #1 \right\rVert_{#2}}

\newcommand{\inprod}[3]{\left( #1, #2 \right)_{#3}}

\newcommand{\abs}[1]{\left\lvert #1 \right\rvert}

\DeclareMathOperator*{\argmax}{arg\,max}
\theoremstyle{definition} 

\newtheorem{assumption}{Assumption}[section]

\numberwithin{equation}{section}

\numberwithin{equation}{section}

\newcommand{\rv}{\right\rangle}

\newcommand{\bR}{\bm{R}}
\newcommand{\bzbar}{\overline{\bm{z}}}
\newcommand{\ba}{\bm{a}}
\newcommand{\bL}{\bm{{\cal L}}}

\newcommand{\bphi}{\bm{\phi}}

\newcommand{\bpsi}{\bm{\psi}}

\newcommand{\bH}{\bm{H}}
\newcommand{\bzero}{\bm{0}}

\newcommand{\bV}{\bm{V}}

\newcommand{\Th}{{\cal T}_h}

\graphicspath{{figures/}}

\def\rv{\xi}         
\def\brv{\bxi}        

\def\banach{L^{\infty}(0,T)\times \bH^1(\Omega)}

\def\normsub{{L^\infty(\dom;\banach)}}

\def\numpts{n}

\newtheorem{theorem}{Theorem}[section]
\newtheorem{lemma}[theorem]{Lemma}

\begin{document}

\begin{frontmatter}
\title{Enhancing adaptive sparse grid approximations and improving refinement strategies using adjoint-based a posteriori error estimates}
\author[sandia]{Jakeman, J.D.\fnref{fn1}\corref{cor1}}\ead{jdjakem@sandia.gov}
\author[sandia]{Wildey, T.\fnref{fn1}}
\address[sandia]{Sandia National Laboratories, 
		Albuquerque, NM 87185, 
		United States}
\cortext[cor1]{Corresponding author}
\fntext[fn1]{Sandia National Laboratories is a multi-program laboratory 
managed and operated by Sandia Corporation, a wholly owned subsidiary of Lockheed 
Martin Corporation, for the U.S. Department of Energy’s National Nuclear Security 
Administration under contract DE-AC04-94AL85000.}

\begin{abstract}
In this paper we present an algorithm for adaptive sparse grid approximations
of quantities of interest computed from discretized partial differential equations.
We use adjoint-based a posteriori error estimates of the physical discretization error and the 
interpolation error in the sparse grid to enhance the sparse grid approximation and to drive adaptivity of the sparse grid. 
Utilizing these error estimates provides significantly more accurate functional values for random samples of 
the sparse grid approximation. We also demonstrate that alternative refinement strategies based upon a posteriori error estimates
can lead to further increases in accuracy in the approximation over traditional hierarchical surplus based strategies. 
Throughout this paper  we also provide and test a framework for balancing the physical discretization error with the stochastic interpolation
error of the enhanced sparse grid approximation. 
\end{abstract}

\end{frontmatter}

Partial Differential Equations (PDE) are used to simulate a wide
range of phenomenon and are often used to inform design decisions and to estimate risk in systems 
with large human and/or financial impact but with limited capacity for experimentation. 
Given the importance of these 
applications the ability to accurately quantify uncertainty in model 
predictions is essential. 

Most uncertainty quantification (UQ) studies focus on estimating parametric 
uncertainty. In such analyses, the uncertainty in the input data, such 
as model coefficients, forcing terms etc, is usually represented through 
a finite number of random variables with a known probability distribution.
The goal of the study is then to compute the effect of the varying input data
on the system response, and in many cases, to calculate the statistics of the 
response.

The accuracy to which uncertainty can be
quantified is limited by the computational resources available to resolve the
governing equations. Many models require vast amounts of computational
effort and thus the number of model evaluations that can be used to interrogate
the uncertainty in the system behavior is limited. Consequently a significant portion 
of methods developed for uncertainty quantification (UQ) in recent years have 
focused on constructing surrogates of expensive simulation models using only a 
limited number of model evaluations. 

The most widely adopted approximation methods are based on
generalized polynomial chaos (PC) expansions~\cite{GhanemSpanos,xiu02a}, 
sparse grid interpolation~\cite{jakeman2013localuq,ma09} and
Gaussian Process (GP) models~\cite{rasmussen2006}. The performance of these methods is problem dependent and in 
practice it is difficult to estimate the accuracy of the approximation constructed.
Cross-validation is one means of
estimating the accuracy of the approximation, however the accuracy of the cross-validation
prediction of the error is limited.  Moreover, cross validation is not readily applied 
for approximation methods which require structured model samples, such as 
sparse grid interpolation and many forms 
of pseudo-spectral projection.

In this paper we utilize sparse grid interpolation to approximate model responses. Sparse grids
can be built using local or global basis functions and have well established and effective 
adaptivity procedures which can be leveraged in conjunction with good error estimates to 
concentrate computational effort to resolving important dimensions and/or regions of the 
random parameter space. Unlike regression based PCE or Gaussian Process models, sparse grids
can be used regardless of the computational budget of the UQ analysis. For example 
sparse grids can be used to approximate a model response using tens to millions of model
runs, whereas the aforementioned alternatives have upper limits in the low thousands imposed
by the need to solve large linear systems.

Throughout this paper, we will use $J(\brv)$ to denote the exact response (quantity of interest) from a partial differential equation that depends
on the unknown variable $\brv$.  When solving PDEs using techniques such as the finite element method the physical discretization error will be non-zero. 
We use $J_h(\brv)$ to denote the response from the discretized model.  As previously mentioned, solving the discretized model is often computationally expensive and therefore we need to consider a surrogate approximation of $J_h(\brv)$, which we denote $J_{h,n}(\brv)$. 
Given these approximations, the error in the response can be decomposed into two components
\begin{equation}\label{eq:error-decomposition}
 \norm{\qoi-J_{h,n}(\brv)}{} \le
 \underbrace{\norm{J(\brv) - J_h(\brv)}{}}_{I} +
 \underbrace{\norm{ J_h(\brv) - J_{h,n}(\brv)}{}}_{II}
\end{equation}
where:
I is the finite element {\it physical discretization} error; and 
II is the {\it stochastic approximation} error introduced by approximating the quantity of interest by a sparse grid
interpolant.

Recently, a posteriori error estimation has arisen as a promising approach
to estimate the error in approximate input-output relationships. Adjoint-based a posteriori error
estimation was originally developed to estimate error in numerical approximations of 
deterministic Partial Differential Equations (PDE)~\cite{BangerthRannacher,eehj_book_96,GilesSuli,oden2001goal}, but recent modifications allow similar ideas to be used to estimate error
in approximations of quantities of interest from PDEs with uncertain parameters.
This relatively new approach, introduced in~\cite{BDW} and further analyzed in~\cite{butler_constantine_wildey1,Butler2013}, 
is based on goal-oriented adjoint-based error estimates and is used to predict error in samples
of a response surface approximation of a specific quantity of interest.  
Similar to standard adjoint-based error estimation procedures, this new approach includes 
the physical discretization error if the adjoint problem is approximated in a higher-order discretization space.
However, the error estimate from this new approach also contains an approximation of the error in
the stochastic discretization due to the evaluation of the response surface model rather than the PDE.
In~\cite{butler_constantine_wildey1,Butler2013}, it was shown that, for spectral and pseudo-spectral Galerkin approximations, this estimate of the stochastic interpolation error is higher-order even if 
a low order approximation of the adjoint is used for the stochastic approximation.  

In general, it is inefficient to reduce the stochastic error to a level below the error introduced by the deterministic discretization. 
Much of the existing literature focuses on minimizing the stochastic approximation error, 
however only a few attempts have been made to discuss or 
account for the combined effect of deterministic and stochastic approximation error.
Error bounds for the stochastic approximation error for isotropic sparse grid approximations of 
elliptic PDEs using Clenshaw-Curtis or Gaussian abscissa 
are given in ~\cite{nobile08a}. In this paper, we use adjoint-based error estimates to ensure
that the error in the stochastic approximation is never significantly less than the physical discretization error.

Our goal in this paper is to utilize adjoint-based a posteriori error estimates to efficiently compute pointwise approximations of specific quantities of interest, 
usually computed from PDE solutions, using adaptive sparse grid approximations.
Specifically, we aim to
\begin{itemize}
 \item Provide theoretical bounds on the error in a posteriori enhanced Clenshaw-Curtis sparse grids.
 \item Numerically demonstrate the enhancement results in \cite{butler_constantine_wildey1} extend to adaptive sparse grid approximations.
 \item Present new refinement strategies for sparse grids based on a posteriori error estimates.
 \item Present a strategy for reducing the cost of computing a posteriori error estimates.
 \item Provide a criterion to stop sparse grid refinement when the stochastic approximation error of the sparse grid is approximately equal 
  to the physical discretization error.
\end{itemize}

The remainder of this paper is organized as follows.  Section~\ref{sec:notation} introduces the general model problem we are interested in.
Sparse grid approximation is reviewed in Section~\ref{sec:sparse-grids} and we recall the standard adjoint-based posteriori error 
analysis for deterministic PDEs in Section~\ref{sec:error-analysis}.
In Section~\ref{sec:SGError} we formulate an a posteriori error estimate for samples of a sparse grid surrogate
and derive theoretical bounds on the error in the a posteriori error estimate.
Section~\ref{sec:sparse-grid-adaptivity}, presents new adaptive strategies for sparse grid refinement that leverage a posteriori error estimates and 
Section~\ref{sec:sg-approx-of-error} introduces the sparse grid approximation of the error estimate and our stopping criteria based on an estimate the physical discretization error.  
Numerical results are presented in Section~\ref{sec:results} and our conclusions are presented in Section~\ref{sec:conclusions}.

\section{General nonlinear problem and notation}\label{sec:notation}

We consider the following system,
\begin{equation}\label{eq:genlin}
\frac{\partial \sol}{\partial t}+\bA(\brv;\sol) = \bzero. 
\end{equation}
In this paper, we consider systems of ordinary differential equations (ODEs) as well as systems of partial differential equations (PDEs).
We focus most of our attention on the numerical approximation of PDEs since the corresponding analysis for the ODE case is a straightforward simplification.
In the PDE case, \eqref{eq:genlin} is defined on $\Omega \times (0,T]$ where $\Omega \subset \mathbb{R}^s$, $s=1,2,3$, is a polygonal
(polyhedral) and bounded domain with boundary $\partial \Omega$.
The random parameter $\brv$ takes values in $\dom \subset \mathbb{R}^d$ and reflects uncertainty in model and source parameters. 
The solution operator's dependency on $\brv$ implies that $\sol:=\sol(\bx,t,\brv)$ is also uncertain and may be modeled 
as a random process for which we will construct a sparse grid approximation.

Before introducing assumptions on the properties of $\bA$ and $\sol$ let us first introduce some notation.
Let $(\cdot,\cdot)_{\cal D}$ denote the inner product of $\bm{L}^2({\cal D})$
and if the domain of integration is clear from the context, we suppress the index ${\cal D}$.
Define $\bV$ be a Sobolev space where for any non-negative integer $m$ we recall 
$$H^{m}(\Omega) = \{v \in L^2(\Omega)\,;\, \partial^k v \in L^2(\Omega)\, \forall\, |k|\le m\},$$
equipped with the following seminorm and norm
$$|v|_{H^{m}(\Omega)} = \left[\sum_{|k|=m} \int_\Omega |\partial^k
v|^2\,dx\right]^{1/2}\ ,\
\|v\|_{H^{m}(\Omega)} = \left[\sum_{0\le |k| \le m}
|v|^2_{H^{k}(\Omega)}\right]^{1/2}.$$
We use $\bH^m(\Omega)$ to denote the obvious generalization to vector valued functions. 

In this paper we assume that $\bA$ is convex and has smooth second derivatives.
Specific examples of $\bA$ and $\sol$ will be given in subsequent sections.  
We assume that sufficient initial and boundary conditions are provided so that \eqref{eq:genlin} is well-posed in the
sense that there exists a solution for a.~e. $\brv \in \dom$. Specifically we assume that (\ref{eq:genlin}) has an equivalent variational formulation seeking $\sol \in \bV$, for any $\brv \in \dom$, 
such that
\begin{equation}\label{eq:genweak}
 \int_0^T \left[ \left(\frac{\partial \sol}{\partial t},\bw\right) + \ba(\brv;\sol,\bw) \right] \ dt = \bzero, \quad \forall \bw \in \bV, 0<t\leq T,
\end{equation}
where $\ba(\brv;\cdot,\cdot)$ is the weak form of $\bA$.

Now let $\Th$ be a conforming partition of $\Omega$, composed of $N_T$ closed convex volumes of maximum diameter $h$.
An element of the partition $\Th$ will be
denoted by $T_i$ where $h_i$ stands for the diameter of $T_i$ for $i=1,2,\ldots,N_T$.
We assume that the mesh is regular in the sense of Ciarlet~\cite{ref:Ciar78} and take $\Th$ to be a conforming finite element mesh consisting of simplices or parallelopipeds.
Denoting the space of continuous, piecewise functions of degree $q$ (simplices) or total degree $q$ (parallelopipeds) over the spatial domain by
\[
	\bV_h^{(q)} = \left\{ \bv \in C(\Omega)\cap \bH^1(\Omega): \forall E\in \Th, \bv|_E \in \mathbb{P}^q(E)\right\}.
\]
we then use a semi-discrete variational formulation to seek $\solh \in \bV_h^{(q)}$, for any $\brv \in \dom$, such that
\begin{equation}\label{eq:genweakdisc}
 \int_0^T \left[ \left(\frac{\partial \solh}{\partial t},\bw\right) + \ba(\brv;\solh,\bw) \right] \ dt = \bzero, \quad \forall \bw \in \bV_h^{(q)}, 0<t\leq T.
\end{equation}

A fully discrete scheme for any $\brv \in \dom$ can be obtained by letting $I_n=(t_{n-1},t_n)$ and time steps $k_n = t_n - t_{n-1}$ 
denote the discretization of $[0,T]$ as $0 = t_0 < t_1 < \cdots < t_N = T$.
In this paper, we do not focus on any particular time integration scheme.  We do assume that the
fully discrete approximation is given as a polynomial in time, e.g., 
\[  \solh \in \bV_h^{(q)} \times \mathcal{W}^{(r)}, \]
where  
\[\mathcal{W}^{(r)} = \left\{w \in C([0,T]): w \in \mathbb{P}^r(I_n), \forall I_n \right\},\]
denotes the space of continuous piecewise polynomial functions of degree $r$.  
To make the notation less cumbersome, 
we generally use $\solh$ to denote the fully discrete solution unless otherwise noted.
The choice of continuous polynomials in time is merely for convenience.  Discontinuous polynomials may also be used
and require a straightforward modification to the a posteriori error analysis involving jump terms
at each time node \cite{eehj_book_96,ES00}.

In cases where a unique and sufficiently regular solution to \eqref{eq:genweak} exists, one can obtain the following error bound (see e.g.~\cite{Brennerscott,Gunzflow,ES00,eehj_book_96} for specific examples) for any $\brv \in \dom$,
\begin{equation}\label{eq:deterrbnd}
 \|\sol - \solh\|_{L^{\infty}((0,T);\bH^1(\Omega))} \leq C(\sol) (h^s + k^{\beta}),
\end{equation}
where $C(\sol)$ depends but does not depend on $k = \max_n k_n$ or $h$.  The parameters, $s$ and $\beta$, are determined by the regularity of the solution and the order of accuracy of the spatial and temporal discretizations.  We emphasize that obtaining error bounds of this sort is nontrivial and may not be possible in many cases.  Throughout the remainder of this paper we assume that \eqref{eq:genweak} does have a unique and sufficiently regular solution to obtain \eqref{eq:deterrbnd}.
\section{Sparse grid approximations}
\label{sec:sparse-grids}
Sparse grid stochastic collocation has been shown to provide efficient and accurate approximation of stochastic 
quantities~\cite{jakeman2013localuq,ma09,nobile08a}. Here we adopt this technique to 
construct approximations of functionals of the discretized solution, $\qoih:=J(\solh(\brv))$.
We typically do not know the closed form of $\qoih$ and only require that we
can evaluate $\qoih$ at arbitrary points in $\dom$ by evaluating the discretized model.
For simplicity we will restrict attention to consider stochastic collocation problems
characterized by variables $\bxi$ with finite support normalized to fit in the domain
$\dom=[0,1]^d$. However the technique proposed here can be applied to semi or unbounded
random variables using the methodology outlined in~\cite{Jakeman10Epistemic}. 

Sparse grids~\cite{bungartz04} approximate $\qoih$ via a 
weighted linear combination of basis functions
\begin{equation}\label{eq:sgInterpolantsum}
\intp{n}{\qoih} := \qoinh = \sum_{k=1}^n v_{k}\, \Psi_{k}(\bxi)
\end{equation}
The approximation is constructed on a set of anisotropic grids $\Xi_\bl$ on
the domain $\dom$ 
where $\bl=(l_1,\ldots,l_d)\in \mathbb{N}^d$ is a multi-index denoting the level
of refinement 
of the grid in each dimension. These rectangular grids are Cartesian product
of nested one-dimensional grid points 
$\Xi_{l}=\{\xi_{l,i}:i<0\le i \le m_{l}\}$
\[ \Xi_\bl = \Xi_{l_1} \times \cdots \Xi_{l_d}
\]
Typically when approximating $\qoih$ with a smooth dependence on $\brv$, $\Xi_{l}$
 are chosen to be the nested Gaussian quadrature rules associated with the distribution of 
$\xi_i$. For example the Gauss-Patterson rule is used for uniform 
variables and Genz-Keister rule for Gaussian variables. For functions of 
lower regularity $\Xi_{l}$ are typically chosen
to be equidistantly spaced. The number of points $m_{l}$ of a one-dimensional 
grid of a given level is dependent on the growth rate of the quadrature rule chosen.

The multivariate basis functions $\Psi_k$  are a tensor product of one
dimensional basis functions. 
Adopting the multi-index notation use above we have
\begin{equation}
\label{eq:multi-dimensional-basis} 
\Psi_{\bl,\bi}(\bxi)=\prod_{n=1}^d \psi_{l_n,i_n}(\xi_n)
\end{equation}
where $\bi$ determines the location of a given grid point. There is a one-to-one 
relationship between $\Psi_k$
in~\eqref{eq:sgInterpolantsum} and $\Psi_{\bl,\bi}$ and each $\Psi_{\bl,\bi}$
is 
uniquely associated with a grid point 
$\bxi_{\bl,\bi}=(\xi_{l_1,i_1},\ldots,\xi_{l_d,i_d})\in \Xi_\bl$. 
Many different one-dimensional basis functions $\psi_{l_n,i_n}(\xi_n)$ 
can be used. If $\qoih$ has a smooth dependence on $\brv$ then the best choice
 is the one-dimensional Lagrange
polynomials. If local approximation is required one can use a multi-dimensional
piecewise polynomial basis~\cite{bungartz04}.

The multi-dimensional basis~\eqref{eq:multi-dimensional-basis} span 
the discrete spaces $V_\bl \subset L^2(\dom)$
\[
 V_\bl=\text{span}\left\{ \Psi_{\bl,\bi}\,:\bi\in\cK_\bl\right\} \quad 
\cK_\bl=\{\bi:i_k=0,\ldots,m_{l_k}\,,k=1,\ldots,d\}
\]
These discrete spaces can be further decomposed into 
hierarchical difference spaces
\[
 W_\bl={V_\bl} \setminus V_{\bl}\, \bigoplus_{n=0}^d V_{\bl-\be_n}
\]
The subspaces $W_\bl$ consists of all basis functions $\Psi_{\bl,\bi} \in V_\bl$
which are not included in any of the spaces $V_\bk$ smaller than $V_l$, 
i.e. with $\bk<\bl$.
These hierarchical difference spaces can be used to decompose the input space
such that
\begin{equation*}
\label{eq:vspace}
V_\bl=\bigoplus_{\bk\le \bl}W_\bl\quad\text{and}\quad
 L^2(\dom)=\bigoplus_{k_1=0}^{\infty}\cdots\bigoplus_{k_d=0}^{\infty}
W_\bk=\bigoplus_{\bk\in\mathbb{R}^d}W_\bk
\end{equation*}

For numerical purposes we must truncate the number of difference spaces used to
construct $V$.  Traditional isotropic 
sparse grids can be obtained by all hierarchical subspaces $W_\bl$ with
an index set that satisfy
\begin{equation}
\label{eq:sg-isotropic-index-set}
 \cL = \{\bl:|\bl|_1\le l\}
\end{equation}

Given a truncation, such as the a-priori one above or one which has been determined adaptively,
$\qoih$ can be approximated by
\begin{equation}\label{eq:sgInterpolant}
\qoinh=\sum_{\bl\in\cL}  J_\bl,\quad J_\bl=
\sum_{\bi\in\cI_\bl} v_{\bl,\bi}\, \Psi_{\bl,\bi}(\bxi)
\end{equation}
where $\cI_\bl = \{\bi:\Psi_{\bl,\bi}\in W_\bl \}$.

Here we note that the $v_{\bl,\bi}$ are the coefficient values of the hierarchical product 
basis, also known as the hierarchical surplus. The surpluses are simply 
the difference between the function value and the sparse grid approximation at
a point, not already in the sparse grid. That is 
\[
v_{\bl,\bi} = J(\bxi_{\bl,\bi}) - J_n(\bxi_{\bl,\bi}),\quad 
\cL \cup \bl = \emptyset
\]

\section{Adjoint-based a posteriori error estimation}\label{sec:error-analysis}


The goal of a simulation is often to accurately estimate a relatively small number of quantities of interest $\qoi$.
Adjoint-based error analysis relates the error ($\be = \sol - \solh$) in a quantity of interest to
a computable weighted residual.  We are usually interested in estimating the error in a numerical 
approximation computed using a discretization of a variational formulation, but it is often 
intuitive to start the discussion with strong form adjoint operators.
We assume that we are interested in estimating the error in a linear functional of the solution,
\begin{equation}\label{eq:linearfunc}
\qoi = \left(\bpsi_T(x),\sol(x,T,\brv)\right) + \int_0^T \left(\bpsi,\sol(x,t,\brv) \right) \ dt, 
\end{equation}
where $\bpsi_T(x)$ is used to compute a linear functional at $t=T$, and $\bpsi = \bpsi(x,t)$
is used to compute time-averages.  Typically, either $\bpsi_T$ or $\bpsi$ are zero, but this is not
required.  Nonlinear functionals and linear functionals at $t\neq T$ can also be used, but this
is omitted here for the sake of simplicity.


For a given $\brv \in \dom$, the linear adjoint operator in strong form is defined via the duality relation
\begin{equation}
 \int_0^T \left[ \left(\frac{\partial \sol}{\partial t},\bv\right) + \left(\bL(\brv) \bv, \bw \right)\right] \ dt = 
 \int_0^T \left[ \left(-\frac{\partial \bv}{\partial t},\sol\right) + \left(\bv, \bL^*(\brv)\bw \right) \right] \ dt 
\label{eqn:adjoint-def}
\end{equation}
where $\bL(\brv)$ is a given linear operator. For a general nonlinear PDE one approach to define the
linear operator $\bL(\brv)$ is to assume $\bA(\brv;\cdot)$ is convex and use the Integral Mean Value Theorem
yielding
\[\bA'(\brv,\bzbar;\be) = \bA(\brv;\sol) - \bA(\brv;\solh)\]
where $\bzbar$ lies on the line connecting $\sol$ and $\solh$, and $\be=\sol-\solh$.
In practice, $\bzbar$ is unknown so we linearize around $\solh$ giving, 
\[\bL(\brv) \be = \bA'(\brv,\solh;\be) = \bA(\brv;\sol) - \bA(\brv;\solh) + \bR(\brv;\be,\solh,\bzbar),\]
where $\bR(\brv;\be,\solh,\bzbar)$ represents the remainder.
Since $\solh - \bzbar \approx \be$, it is common to assume that the remainder is a higher 
order perturbation term and can be neglected \cite{BangerthRannacher,Beck_Ran_opt_a_post_FE_01,ES00}.
Notice that the operator $\bL(\brv)$ is often the same linear operator used in computing the step in Newton's method. 
This fact is often exploited to ease construction of the discrete adjoint operator.

The error in a linear functional can be represented using the definition of the adjoint:
\begin{multline} \label{eq:var-err-rep}
\qoi - \qoih = \left(\bphi(x,0,\brv),\sol(x,0,\brv) - \solh(x,0,\brv)\right) \\
 -\int_0^T \left(\left(\frac{\partial \solh}{\partial t},\bphi \right)+\ba(\brv;\solh,\bphi)\right) \ dt 
 + \int_0^T \bR(\brv;\be,\solh,\bzbar,\bphi) \ dt,
\end{multline}
where $\bphi:=\bphi(x,t,\brv)$ is defined by the adjoint problem
\begin{align}\label{eq:adjoint-strong}
-\frac{\partial \bphi}{\partial t}+\bL^*(\brv) \bphi &= \bpsi(x,t), \\
                                         \phi(x,T) &= \bpsi_T(x).
\end{align}
If the adjoint solution, $\bphi$, is given, then the error representation in Eq.~\eqref{eq:var-err-rep} is easily evaluated if we neglect
the higher order remainder term, see \cite{BangerthRannacher,Beck_Ran_opt_a_post_FE_01,ES00,Estep:2010:PEE:1752623.1752844} for a complete discussion of this remainder term.
However, the solution to the adjoint problem Eq.~\eqref{eq:adjoint-strong} is usually not given explicitly 
and we must approximate the solution using an appropriate discretization.  In this paper, we approximate the adjoint solution using a finite element method.  Since we are also interested in estimating the physical discretization error, 
we use a higher-order approximation in the spatial domain. 
This gives a fully computable error representation:
\begin{multline} \label{eq:error_rep}
\varepsilon(\solh,\bphi_h) = \left(\bphi_h(x,0,\brv),\sol(x,0,\brv) - \solh(x,0,\brv)\right) \\
 -\int_0^T \left(\left(\frac{\partial \solh}{\partial t},\bphi_h \right)+\ba(\brv;\solh,\bphi_h)\right) \ dt.
\end{multline}

In general, the definition of the adjoint problem also requires appropriate boundary conditions and the error representation \eqref{eq:var-err-rep} also includes boundary terms to account for errors
made in approximating the boundary conditions on the forward problem.
We omit these terms in this paper for the sake of simplicity and refer the interested reader to the standard references, e.g.~\cite{BangerthRannacher,eehj_book_96,GilesSuli,oden2001goal}, for a complete discussion of these additional terms. 

\section{Adjoint-based error estimates for samples of a sparse grid surrogate}\label{sec:SGError}

In many practical situations, we are interested in computing probabilistic quantities, such as the probability of a particular event, from the 
surrogate approximation. Computing probabilistic quantities usually requires sampling the surrogate according to the distribution of the random parameters over $\dom$.  In order to have confidence in our estimates of these probabilistic quantities, 
we must consider the error in each sample of the surrogate.  

In this paper, we employ the technique introduced in \cite{BDW} and further analyzed in \cite{butler_constantine_wildey1,Butler2013} to estimate the error in a sample of the sparse grid surrogate of the quantity of interest.  


Given sparse grid approximations of the forward $\solhn$ and adjoint $\bphi_{h,n}$ solutions we can compute the following approximate
error estimate
\begin{equation}\label{eq:erroreval}
\qoi - \qoihn \approx \eest , 
\end{equation}
where
\begin{multline*}
 \eest = \left(\bphi_{h,n}(x,0,\brv),\sol(x,0,\brv) - \solhn(x,0,\brv)\right) \\
- \int_0^T \left(\left(\frac{\partial \solhn(x,t,\brv)}{\partial t},\bphi_{h,n}(x,t,\brv) \right)+\ba(\brv;\solhn(x,t,\brv),\bphi_{h,n}(x,t,\brv))\right) \ dt. 
\end{multline*}


The error in $\eest$ is dependent on the errors in the approximation $\solh$ and $\bphi_h$ which are in turn
dependent on the regularity of $\sol$ and $\bphi$ with respect to the random variables $\brv$. 
For the analysis presented in this paper we will make the following regularity assumption.
\begin{assumption}[Regularity Assumption]\label{assumption:regularity}
Let $\Gamma_{\rv}^1$ be the one dimension range of a random variable $\xi$ and 
$$\Gamma_{\brv}=\prod_{i=1}^d\Gamma^1_{\xi_i},\quad\Gamma^{(-n)}_{\brv}=\prod_{\substack{i=1\\i\neq n}}^d\Gamma^1_{\xi_i}$$
For each $\rv_n\in\Gamma_n^1$ there exists $\tau_n>0$ such that the function $\solh(\rv_n,\brv^{(-n)},\bx):\Gamma^1_{\rv_n}\rightarrow C^0(\Gamma^{(-n)}_{\brv};\banach)$ as a function of $\rv_n$, 
admits an analytical extension $\sol(y,\brv^{(-n)},\bx)$, $y\in\mathbb{C}$ in the region of the complex plain
\[
 \Sigma(\Gamma^1_{\rv_n};\tau_n)\equiv\{y\in\mathbb{C},\text{dist}(y,\Gamma^1_{\rv_n})\le\tau_n\},
\]
where $\brv^{(-n)}\in\Gamma^{(-n)}_{\brv}$. Moreover $\forall y\in\Sigma(\Gamma^1_{\rv_n};\tau_n)$
\[
\norm{\solh(y)}{C^0(\Gamma^{(-n)}_{\brv};\banach)}\le\lambda
\]
for a constant $\lambda$ independent of $n$.
\end{assumption}
This assumption is valid for a number of problems. Indeed the 
the linear elliptic diffusion equation analyzed in Section~\ref{sec:SGError} satisfies assumption~\ref{assumption:regularity}~\cite{babuska2007}.
Using this assumption we can obtain the following result.
\begin{lemma}[\cite{nobile2008}]
\label{lemma:sg-point-wise-error}
For functions $\sol \in C^0(\dom;\banach)$ satisfying Assumption~\ref{assumption:regularity} 
the isotropic Smolyak formula (2.8) based on Clenshaw-Curtis abscissas, and with $n$ points, satisfies:
\begin{align*}
 \norm {\solh-\solhn}{\normsub} & \le C_1(\sigma_1) \numpts^{-\mu_1} \\
 \mu_1 &:=\frac{\sigma_1}{1+\log{2d}}
\end{align*}
where the constant $C_1(\sigma)$ depends on the size of the analyticity region $\sigma_1$ but not on the number of points in the sparse grid.
\end{lemma}

Combining Lemma~\ref{lemma:sg-point-wise-error} and taking the maximum of \eqref{eq:deterrbnd} over $\dom$ yields the following result.
\begin{lemma}[\cite{nobile2008}]
\label{lemma:total-point-wise-error}
For functions $\sol \in C^0(\dom;\banach)$ satisfying Assumption~\ref{assumption:regularity}, 
 isotropic Smolyak formula (2.8) based on Clenshaw-Curtis abscissas satisfies:
\begin{align*}
 \norm {\sol-\solhn}{\normsub} & \le C_1(\sigma_1) \numpts^{-\mu_1} + C_1(\sol)\left(h^{s_1}+k^{\beta_1}\right)
\end{align*}
where the constant $C_1(\sigma_1)$ depends on the size of the analyticity region $\sigma_1$ but not on the number of points in the sparse grid and $C_1(\sol)$
depends on the solution $\sol$ but not $h$ or $k$.
\end{lemma}

We can also apply Lemmas~\ref{lemma:sg-point-wise-error} and \ref{lemma:total-point-wise-error} to the adjoint solution, $\bphi$, and the numerical approximation of the adjoint solution, $\bphi_{h,n}$, giving
\begin{align}\label{eq:adj_error_bound}
 \norm {\bphi-\bphi_{h,n}}{\normsub} & \le C_2(\sigma_2) \numpts^{-\mu_2} + C_2(\bphi)\left(h^{s_2}+k^{\beta_2}\right).
\end{align}
Note that we account for the fact that the constants and convergence parameters may be different for the adjoint solution.

Next, we prove that the error estimate, $\eest$, converges at a faster rate than $\qoihn$.  First, we note that we can approximate
\[\qoi - \qoihn \approx \varepsilon(\solhn,\bphi).\]
This approximation is exact for linear problems.  For nonlinear problems, the error in this approximation is a high-order linearization term
and is usually neglected.  Thus, we prove that $\eest$ is a high-order approximation of $\varepsilon(\solhn,\bphi)$.
\begin{theorem}
\label{thm:error_bound}
For functions $\sol, \bphi \in C^0(\dom;\banach)$ satisfying Assumption~\ref{assumption:regularity} and using 
 isotropic Smolyak formula (2.8) based on Clenshaw-Curtis abscissas, the error estimate satisfies:
\begin{align*}
 \|\varepsilon(\solhn,\bphi) - \eest\|_{L^{\infty}(\dom)} & \le \left(C_1(\sigma_1) \numpts^{-\mu_1} + C_1(\sol)\left(h^{s_1}+k^{\beta_1}\right)\right) \\
& \times \left(C_2(\sigma_2) \numpts^{-\mu_2} + C_2(\bphi)\left(h^{s_2}+k^{\beta_2}\right)\right) 
\end{align*}
where the constants $C(\sigma)$ and $C_2(\sigma_2)$ depend on the size of the analyticity regions $\sigma_1$ and $\sigma_2$ respectively, but not on the number of points in the sparse grid.
In addition, $C_1(\sol)$ and $C_2(\bphi)$
depend on the solutions $\sol$ and $\bphi$ respectively, but not $h$ or $k$.
\end{theorem}

\begin{proof}
 We use the fact that $\varepsilon(\cdot,\cdot)$ is linear in the second argument to write,
\[\left|\varepsilon(\solhn,\bphi) - \eest\right| = \left|\varepsilon(\solhn,\bphi - \bphi_{h,n})\right|.\]
It is straightforward (see e.g.~\cite{Butler2013,oden2001goal,BangerthRannacher}) to show that
\begin{align*}
\left|\varepsilon(\solhn,\bphi - \bphi_{h,n})\right| \leq & \ \norm {\sol-\solhn}{\normsub} \\
& \times \norm {\bphi-\bphi_{h,n}}{\normsub}.
\end{align*}
A direct application of Lemma~\ref{lemma:total-point-wise-error} and \eqref{eq:adj_error_bound} proves the result.
\end{proof}

We are often interested in using $\eest$ to improve the approximation of the quantity of interest.  To this end, we define
the enhanced quantity of interest to be
\[\qoien(\brv):=\qoihn + \eest.\]
Theorem~\ref{thm:error_bound} implies that $\qoien(\brv)$ is higher-order approximation to $\qoi$ than $\qoihn$. 
\section{Sparse grid adaptivity}\label{sec:sparse-grid-adaptivity}
Although a vast improvement of full tensor grid approximations the number
of sparse grid points in the isotropic sparse grids, presented in Section~\ref{sec:sparse-grids} and analyzed in Section~\ref{sec:SGError},
 still grows quickly with dimension. To enhance the convergence of sparse grids two types of adaptivity
have been developed. The first type of adaptivity refines the grid 
dimension by dimension~\cite{gerstner03}, greedily choosing points in dimensions that are deemed
by the algorithm to be more important. The second type of adaptivity
refines the sparse grid locally in regions identified as important~\cite{griebel98}.
When using localized basis functions adaptivity can be performed combining the 
strengths of both dimension and regional adaptivity~\cite{jakeman2013localuq}.

\subsection{Dimension adaptivity}
The dimension adaptive algorithm begins with a low-level isotropic sparse grid
approximation with a set of subspaces $\cL$ and active subspaces $\cA$. 
Often $\cL=W_{\bzero}$ and $\cA=\{W_{\be_k},k=1\ldots,d\}$.
We then choose $W_\bl \in \cA$ with the largest error indicator $\gamma_\bl$ and refine
that subspace. The subspace is refined by adding all indices $W_\bk$ with 
$\bk = \bl+\be_n$, $n=1,\ldots,d$ that satisfy the following admissibility criterion
\begin{equation}
\label{eq:gsg_admissibility}
\bl-\mathbf{e}_k\in\cL\text{ for }1\le k\le d,\, l_k > 1
\end{equation}
Each subspace that satisfies~\eqref{eq:gsg_admissibility} is then placed
in the active set $\cA$. This
process continues until a computational budget limiting the number of model samples
(grid points) is reached or a global error indicator drops below a predefined
threshold. Pseudo-code for the dimension adaptive algorithm is shown in Algorithm~\ref{alg:dim-adaptivity}. 

These \texttt{INDICATOR} and \texttt{GLOBAL\_INDICATOR} routines control which subspaces are added
to the sparse grid. The key to adaptive sparse grids working well is greatly influenced by these two 
routines, which respectively provide
accurate estimates of the contribution of a subspace to reducing the error in the
interpolant, and the error in the entire interpolant.

Typically indicator functions are functions of the hierarchical surplus values $v_{\bl,\bi}$ of points in the grid.
Throughout this paper we will use the indicator
\begin{equation}
\label{eq:dim-surplus-indicator}
\gamma_\bl = \sum_{\bi\in\cI_\bl} \abs{ v_{\bl,\bi} }w_{\bl,\bi}
\end{equation}
as the baseline for comparison for the proposed method.
Here $w_{\bl,\bi}$ is the quadrature weight of the grid point $\rv_{\bl,\bi}$. In addition, we use the global indicator
\[
 \eta = \sum_{\bl\in\cA} \gamma_\bl
\]
Hierarchical surplus indicators such as~\eqref{eq:dim-surplus-indicator} require
the model to be evaluated at the associated grid point $\gridpt$ before the surplus and thus indicator can be computed. 
This procedure adds a grid point and then checks if that refinement should have been performed. Efficiency can be gained
by using an error indicator that flags the need for refinement without the need to evaluate the simulation model. 
The a posteriori error estimate~\eqref{eq:erroreval} is one such indicator.

We construct the following error indicators using the a posteriori error estimate for dimension adaptive sparse grids
\begin{equation}
\label{eq:dim-surplus-indicator-alt}
\gamma_\bl^{\varepsilon} = \sum_{\bi\in\cI_\bl} \abs{ \pteest }w_{\bl,\bi},\quad  \eta^{\varepsilon} = \sum_{\bl\in\cA} \gamma_\bl^{\varepsilon}
\end{equation}
The use of these indicators requires a minor modification to the interpolate algorithms presented in the literature~\cite{gerstner03,jakeman2013localuq}.
The modification, present in Algorithm~\ref{alg:dim-adaptivity}, controls when the model is evaluated at a sparse grid point and when that point is added to the 
sparse grid. For surplus based refinement the model must be evaluated at the sparse grid when \texttt{INDICATOR} is called on a subspace.
When using a posteriori refinement the model does not need to be evaluated
until \texttt{REFINE} is called on the point at which time the point is also added to the sparse grid. 
This modification means that when using a posteriori refinement each step in the algorithm only requires evaluation of the target function $f(\brv_{\bl,\bi})\; \forall \brv_{\bl,\bi}\in W$.
In comparison surplus refinement requires evaluation of $f(\brv_{\bl,\bi})\; \forall \brv_{\bl,\bi}\in \cJ$ at each step. Typically the number of points in $W$ will be much less than $\cJ$ thus
a posteriori refinement should be more efficient than surplus based refinement.

\begin{algorithm}
\caption{\texttt{INTERP}[$f(\brv)$,$\cL$,$\cA$,$\tau$,$n$]$\rightarrow f_{n},\cA,\cL$}
For a given $\cL$ the points in the sparse grid are $\Xi:= \bigcup_{\bl\in\cL} \Xi_\bl$.\\
The number of sparse grid points are $N=\#\Xi$\\\\
WHILE NOT \texttt{TERMINATE}[$\cA$,$N$,$\tau$,$n$]
\label{alg:dim-adaptivity}
\begin{itemize}
 \item $W:= \argmax_{W_\bl\in\cA} \gamma_\bl$ {\footnotesize\% Determine the subspace with the highest priority}
 \item $\cA:=\cA\setminus W$ {\footnotesize\% Remove $W$ from the active set}
 \item IF ( NOT using surplus refinement ) Evaluate $f(\brv_{\bl,\bi})\; \forall \brv_{\bl,\bi}\in W$
 \item $\cJ:=$ \texttt{REFINE}[$W$,$\cL$] {\footnotesize\% Find all admissible forward neighbors of $W$}
 \item IF ( using surplus refinement ) Evaluate $f(\brv_{\bl,\bi})\; \forall \brv_{\bl,\bi}\in \cJ$
 \item $\gamma_{\bl}:= $ \texttt{INDICATOR}[$W_\bl$]$\;\forall\; W_\bl\in\cJ$    {\footnotesize\% Calculate the priority of the neighbors}
 \item $\cA:= \cA \cup \cJ$     {\footnotesize\% Add the forward neighbors to the active index set}
\end{itemize}
\end{algorithm}

\begin{algorithm}
\caption{\texttt{TERMINATE}[$\cA$,$N$,$\tau$,$n$]}
\begin{itemize}
 \item $\eta:=$ \texttt{GLOBAL\_INDICATOR}[$\cA$]
 \item IF $\cA = \emptyset$ or $N\ge n$ or $\eta<\tau$ RETURN TRUE
 \item ELSE RETURN FALSE
\end{itemize}
\label{alg:terminate}
\end{algorithm}

\subsection{Local adaptivity}
Locally-adapted sparse grids can be adapted using a conceptually similar method
to that employed for dimension adapted grids. Instead of refining the grid 
subspace by subspace, locally-adaptive grids are adapted point by point. 
When a new point is chosen for refinement it is added to the sparse grid. 
Children of this point are then found and added to an active set $\cA_\text{local}$. 
In this paper we will consider the typical local refinement strategy~\cite{ma09} and 
the simultaneous dimension and local refinement proposed in~\cite{jakeman2013localuq} termed generalized local refinement.
In this manuscript we will refer to the former as traditional local refinement and
the later as generalized local refinement.
Error indicators, are then computed at each
point in $\cA_\text{local}$. The next point chosen for refinement is simply the 
point $\bxi_{\bl,\bi}$ in  $\cA_\text{local}$ with the largest error indicator
$\gamma_{\bl,\bi}$. The local adaptive sparse grid procedure can be obtained 
from Algorithm~\ref{alg:dim-adaptivity} by defining $\cA$ and $\cL$ to contain points
rather than subspaces and $W$ and $W_\bl$ to be grid points not subspaces. The 
\texttt{REFINE} routine must also be changed as discussed in~\cite{jakeman2013localuq}.

Throughout this paper we will use the indicators
\begin{equation}
\label{eq:local-surplus-indicator}
 \gamma_{\bl,\bi} = \abs{ v_{\bl,\bi} }w_{\bl,\bi},\quad \eta = \sum_{\{\bl,\bi\}\in\cA} \gamma_{\bl,\bi}
\end{equation}
as the baseline for comparison for the proposed method.
Again we can use a posteriori error estimates to construct the alternative refinement metric
\begin{equation}
\label{eq:local-eest-indicator}
 \gamma_{\bl,\bi}^{\varepsilon} = \abs{ \pteest }w_{\bl,\bi},\quad  \eta^{\varepsilon} = \sum_{\{\bl,\bi\}\in\cA} \gamma_{\bl,\bi}^{\varepsilon}
\end{equation}
We remark that when using local refinement one should ensure all the ancestors of a point in the sparse grid have previously been added to the sparse grid.  
Missing ancestors must be added recursively before a child is added. 
The use of a posteriori error guided refinement affects when the \texttt{INTERP} algorithm adds ancestors to the sparse grid.  
When employing a posteriori based refinement, function
values at ancestor points are not necessary to compute the error estimate at a child point.  One has the choice to add missing ancestors either when the 
error estimate of a child is computed or later when the child point is added to the sparse grid. We found during our investigation that adding ancestors before 
computing the error estimate produced a more accurate error estimate at the candidate points 
which resulted in a more accurate sparse grid for a fixed number of sample functions. These results are not shown as the authors believed 
this minor point would distract from the main conclusions that can be drawn from the results section.

Throughout the remainder of this manuscript we will use the following notation to denote the various types of refinement strategies. 
Let $\gamma_\bl$ and $\gamma_\bl^{\varepsilon}$ respectively denote, hierarchical surplus and a posteriori based dimension based refinement. Let
hierarchical surplus based traditional local refinement~\cite{ma09} and generalized local refinement~\cite{jakeman2013localuq} be denoted by 
$\gamma_{\bl,\bi}^{\mathrm{trad}}$ and $\gamma_{\bl,\bi}^{\mathrm{gen}}$, respectively. 
Finally we define $\gamma_{\bl,\bi}^{\varepsilon,\mathrm{gen}}$ to be a posteriori error based local generalized refinement.

\section{Sparse grid approximations of the error}\label{sec:sg-approx-of-error}
As previously mentioned, the evaluation of the a posteriori error estimate using~\eqref{eq:erroreval} requires
the full forward and adjoint approximations at $\brv \in \dom$ as well as the evaluation of the space-time
weighted-residual. While this is usually much cheaper than solving the forward problem at $\brv$, the computation
cost in producing the error estimate at $\brv$ should not be neglected.  
If we only need to compute a small number of samples, then this is probably not a significant issue.
However, we often require a very large number of samples of the surrogate,
and, in turn, a very large number of error estimates.  To mitigate this issue, we propose constructing a surrogate of the error estimate. 
This can significantly reduce the number of error estimates required and the surrogate of the error estimate can be chosen 
to have similar accuracy as the point-wise evaluation of \eqref{eq:erroreval}.

Given a sparse grid approximation of the quantity of interest $\qoihn$ we can use a finite number of error estimates $\pteest$ 
at new grid points to build an enhanced sparse grid. Once built the enhanced sparse grid interpolant can be sampled directly during post-processing, 
thus removing any further need to evaluate~\eqref{eq:erroreval}.

The procedure needed to construct an enhanced sparse grid is outlined in Algorithm~\ref{alg:enhanced-interpolate}.
The algorithm consist of two steps. The first step is used to construct sparse grid approximations of $\qoih$,
forward solution $\solh$ and the adjoint solution $\bphi_h$. The second step uses these sparse grids as the basis of a sparse grid approximation of 
$\qoien(\brv)$, Specifically the final $\cA$ and $\cL$ of $\qoihn$ are used as the initial starting point when building $\qoienm$ and are then 
evolved by sampling $\qoien(\brv)$ at new sparse grid points. We emphasize that the second step no longer involves full PDE solves, but only evaluation of the sparse grids approximation
of the PDE solutions and a weighted residual needed to compute $\eest$.

\begin{algorithm}
\caption{ENHANCED\_INTERP[$\qoih$,$\solh(\brv)$,$\bphi_h(\brv)$,$\cL$,$\cA$,$\tau$,$n$]$\rightarrow \qoien$}
Define $C$ to be the ratio of the cost of evaluating the error estimate 
to the combined cost of a forward and adjoint equation evaluation
\begin{itemize}
\item $\qoih_{h,n},\cA_{n},\cL_{n}$ = INTERP[$\qoih$,$\cA,\cL$,$\tau$,$n/4$]
\item $\sol_{h,n}(\brv,\bx,\bt),\cA_{n},\cL_{n}$ = INTERP[$\solh(\brv,\bx,\bt)$,$\cA,\cL$,$\tau$,$n/4$]
\item $\bphi_{h,n}(\brv,\bx,\bt),\cA_{n},\cL_{n}$ = INTERP[$\bphi_h(\brv,\bx,\bt)$,$\cA,\cL$,$\tau$,$n/4$]
\item $\delta_{\max} = \max_{\gridpt\in\cA}\ptdest$
\item $\gamma_{\max} = \max_{\gridpt\in\cA}\gamma_{\bl,\bi}$
\item $\tau_{\varepsilon}=\max( \delta_{\max} , \gamma_{\max}^2 )$
\item $\qoienm(\brv),\cA_{n,m},\cL_{n,m}$ = INTERP[$\qoihn+\eest$,$\cA_{n},\cL_{n}$,$\tau_{\varepsilon}$,$Cn/2$]
\end{itemize}
\% $\sol_{h,n}$ and $\bphi_{h,n}$ are needed to compute $\eest$.
\label{alg:enhanced-interpolate}
\end{algorithm}

\begin{figure}[ht]
\centering
\includegraphics[width=0.45\textwidth]{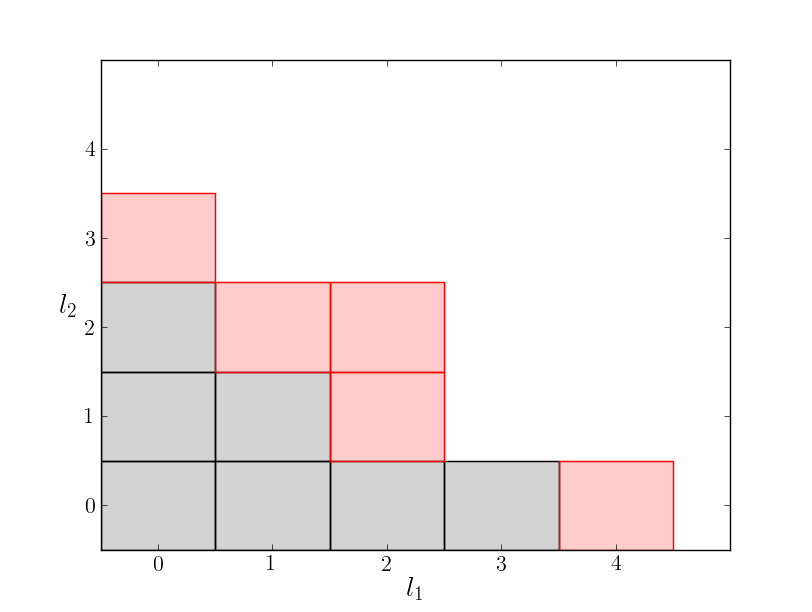}
\includegraphics[width=0.45\textwidth]{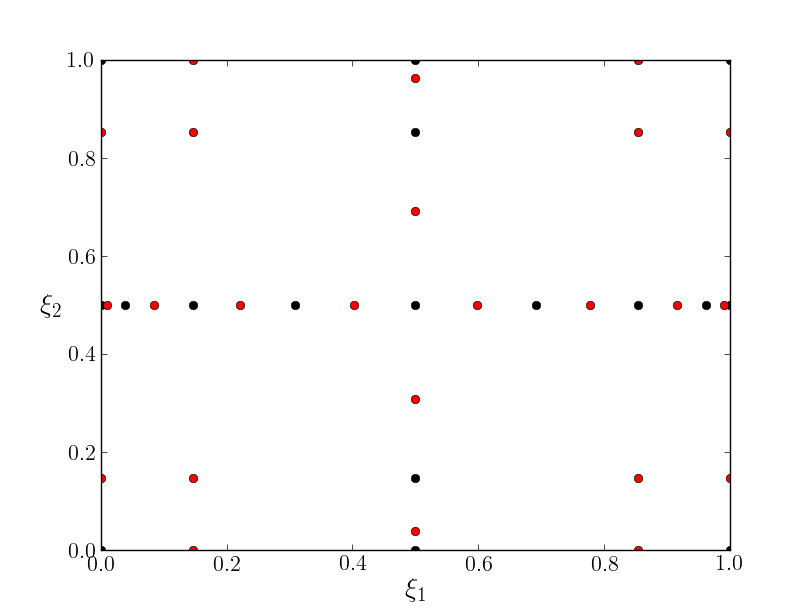}
\caption{The two phases undertaken when building an enhanced dimension adaptive sparse grid interpolant $\qoien$. The first phase consist of building the un-enhanced approximation $\qoihn$. 
The gray boxes on the left represent the 2D sparse grid subspaces $\bl=(l_1,l_2)$ used to build $\qoihn$ and the black points on the right represent the corresponding points 
in that sparse grid. At each of the black points both the forward and adjoint equations are solved. The second phase involves continuing refinement until the sparse grid 
adequately resolves the enhanced function $\qoien$. The red boxes represent the additional subspaces used to build the enhanced approximation $\qoienm$. The red points are 
the associated additional points that are added to the sparse grid. Unlike the black points the red points only require evaluation of the residual (to compute the a posteriori error estimate) and 
the cheaply evaluated approximation $\qoihn$.
The red boxes are not to be confused with active indices which cannot have the structure shown in this figure.}
\label{fig:two-phase-interpolation}
\end{figure}
%
%

\subsection{Adaptivity}\label{subsec:error-adaptivity}

When using dimension-adaptive sparse grids
we begin with an initial grid with subspaces in $\cL_n$ and a set of candidate subspaces $\cA_n$ generated when building the sparse grid approximation 
of the quantity of interest, forward solution and the adjoint solution. When using surplus refinement $\cA_n\subset\cL_n$ and when using 
a posteriori error refinement $\cL_n\cap\cA_n=\emptyset$. This difference arises from the fact that surplus refinement requires the model to be evaluated at all the points in the candidate subspaces
where as a posteriori refinement does not. The set $\cL_n$ at this stage of the algorithm are the gray boxes in Figure~\ref{fig:two-phase-interpolation}.

Given the sets $\cA_n$ and $\cL_n$ we simply use \texttt{INTERP}[$\qoien(\brv)$,$\cL_n$,$\cA_n$,$\tau_{\varepsilon}$,$m$] to construct the enhanced sparse grid $\qoienm:=I_m[\qoien(\brv)]$.
When the algorithm terminates $\cL$ will include the red and gray boxes in Figure~\ref{fig:two-phase-interpolation} and $\cA=\emptyset$.

Aside from noting this minor distinction, constructing $\qoienm$ requires one additional step. 
Before \texttt{INTERP}[$\qoien(\brv)$,$\cL_n$,$\cA_n$,$\tau_{\varepsilon}$,$m$] is called the function value at all points in $\cL_n$ must be updated to include the physical discretization error estimate $\dest$. That is for all 
$\gridpt\in\cL_n$ set the function value at that point to $\ptqoih+\ptdest$. The hierarchical surpluses at each of these points must also be updated accordingly. We do this 
because any new points to the enhanced sparse grid will have function values $\ptqoihn + \pteest$, where $\pteest$ includes the deterministic error estimate. We again remark
that in order to compute the physical discretization error estimate, the physical discretization used to compute the adjoint solution must be higher-order than the method used to compute the forward solution.  This is standard practice in adjoint-based error estimation for deterministic problems.

The method used to construct an enhanced sparse grid using local refinement is similar to that employed when using dimension refinement. Again we must simply define the initial sets
$\cL_n$ and $\cA_n$ and then call \texttt{INTERP} after adjusting the function values and hierarchical surpluses for all $\gridpt\in\cL_n$. 

\subsection{Balancing the stochastic and physical discretization errors}\label{sec:error-balancing}
When quantifying the uncertainty of model predictions, both the deterministic and 
stochastic discretization errors must be accounted for.  It is inefficient to 
reduce the stochastic error to a level below the error introduced by the physical discretization. 
In the following, we assume
that physical discretization, and thus the error induced by the physical discretization, is fixed. 
In order to balance the stochastic error with  this fixed physical discretization error, we must be able to handle two regimes:
first that our computational budget is large enough to drive the stochastic error below the physical discretization error; and secondly the computational budget is insufficient
to force the stochastic error below the physical discretization error. 

In the first regime we must be able to identify when the sparse grid is sufficiently refined.
We use $\tau_{\varepsilon}$ to terminate refinement of $\qoienm$ when the interpolation error is approximately equal to the physical discretization error.
Specifically, we set $\tau_{\varepsilon}$ to be the maximum of the approximate physical discretization error 
$\delta_{\max}$ and the approximate stochastic error $\gamma_{\max}$ in the enhanced sparse grid. 
The error in the error estimate is the product of the error in the approximation of the forward solution and the error in the approximation of the adjoint solution \cite{BDW,butler_constantine_wildey1,Butler2013}.
Consequently a reasonable approximation of the potential accuracy of the enhanced sparse grid is the square of the maximum indicator $\gamma_{\bl,\bi}$ of all points in $\cA$. 
That is we set $\tau_{\varepsilon}=\eta^2$

In the second regime, the goal is to reduce the total error as much as possible within a computational budget of $n$ computational units.  
To do this we must consider the costs of solving the forward problem, solving the adjoint problem, and evaluating the error estimate \eqref{eq:erroreval}.
From our experiments we found it to be advantageous to limit the number of forward 
and adjoint solves to be $1/2$ of the computational budget, and to utilize the remaining effort to evaluate $\eest$ to build the enhanced sparse grid.  
Here we have made the reasonable assumption
that the cost of one forward solve and one adjoint solve are approximately equal, thus the choice of $n/4$ when constructing $\sol_{h,n/4}$ and $\bphi_{h,n/4}$.
The cost in solving the adjoint problem may actually be smaller than the cost in solving the forward problem, even if a higher order method is used for the adjoint, since the adjoint problem is always linear.  

We remark that these termination conditions removes the need for knowledge of what computational budget regime is active.

\section{Results}\label{sec:results}
In this section, we provide several numerical examples to illustrate the proposed methodology. 
In all of the figures presented, the expression left of the colon in the legend denotes the quantity being approximated 
and the expression on the right denotes the type of refinement employed. 
For definitions of the refinement types refer to the end of Section~\ref{sec:sparse-grid-adaptivity}.

\subsection{Heterogeneous diffusion equation}\label{sec:heat-equation}
In this section, we investigate the performance of the proposed methodology when applied to 
a Poisson equation with an uncertain heterogeneous diffusion coefficient. 
Attention is restricted to the one-dimensional physical space to avoid unnecessary 
complexity. The procedure described here can easily be extended to higher physical dimensions.

Consider the following problem with $d \ge 1$ random dimensions:
\begin{equation}\label{eq:hetrogeneous-diffusion}
-\frac{d}{dx}\left[a(x,\brv)\frac{d\sol}{dx}(x,\brv)\right] = 10,\quad 
(x,\brv)\in(0,1)\times I_{\brv}
\end{equation}
subject to the physical boundary conditions
\begin{equation}
\sol(0,\brv)=0,\quad \sol(1,\brv)=0.
\end{equation}
Furthermore assume that the log random diffusivity satisfies
\begin{equation}\label{eq:diffusivityZ}
a(x,\brv)=\bar{a}+\sigma_a\sum_{k=1}^d\sqrt{\lambda_k}\phi_k(x)\rv_k
\end{equation} 
where $\{\lambda_k\}_{k=1}^d$ and $\{\phi_k(x)\}_{k=1}^d$ are, respectively, 
the eigenvalues and eigenfunctions of the covariance kernel 
\begin{equation}\label{eq:heat-eq-qoi}
 C_a(x_1,x_2) = \exp\left[-\frac{(x_1-x_2)^2}{2l_c}\right].
\end{equation}
The variability of the diffusivity field~\eqref{eq:diffusivityZ} is 
controlled by $\sigma_a$ and the correlation length
$l_c$ which determines the decay of the eigenvalues $\lambda_k$. 
In the following we set $d=25$, $\sigma_a=1.$, $l_c=0.1$, $\bar{a} = 0$
and define $\rv_k\in[-1,1]$, $k=1,\ldots,d$ to be independent and identically distributed uniformly random variables.

We are interested in quantifying the uncertainty in the solution $\sol$ at $x=0.5$ caused by the aforementioned random
diffusivity field~\eqref{eq:diffusivityZ}. Accordingly, we define our quantity of interest to be the linear functional 
\[\qoi = \inprod{\psi(x)}{\sol(x,\brv)}{},\quad\psi(x) = C_s\exp(-100(x-.5)^2),\]
where $C_s$ is a scaling constant chosen so that $\int_0^1 \psi(x) \ dx = 1$.
The adjoint problem corresponding to this quantity of interest is given by
\begin{equation}\label{eq:adjhetrogeneous-diffusion}
-\frac{d}{dx}\left[a(x,\brv)\frac{d\bphi}{dx}(x,\brv)\right] = \psi(x),\quad 
(x,\brv)\in(0,1)\times I_{\brv}
\end{equation}
subject to the physical boundary conditions
\begin{equation}
\bphi(0,\brv)=0,\quad \bphi(1,\brv)=0.
\end{equation}

We approximate the solution to~\eqref{eq:hetrogeneous-diffusion} using a standard continuous Galerkin Finite Element Method (FEM)
on a uniform mesh with 100 elements ($h=0.01$) and piecewise linear basis functions. Similarly, to solve \eqref{eq:adjhetrogeneous-diffusion}
we use a continuous Galerkin method with piecewise quadratic basis functions. A higher-order discretization 
is used to solve the adjoint problem so that estimates of the physical discretization error (due to the finite element approximation) can be obtained.

\subsubsection{Adjoint based estimates of the stochastic error}
Here, we demonstrate that sparse grid estimates of the linear functional $\qoihn$ can be made significantly more accurate by
combining these estimates with a posteriori error estimates $\eest$.

Figure~\ref{fig:heat-eq-convergence-d-25-unenhanced-vs-enahnced-brute-only}
compares the errors $\norm{\qoih-\qoihn}{\ell_2(\dom)}$ of un-enhanced sparse surplus-based, dimension-adaptive, Clenshaw-Curtis 
grids with the errors of sparse grid approximations enhanced with a posterior error estimates $\norm{\qoih-(\qoihn+\eest)}{\ell_2(\dom)}$. 
Here and throughout the paper discrete $\ell_2$ errors are computed using 100,000 Latin hypercube samples of $\dom$.
We observe that the rate of convergence
of the un-enhanced sparse grid is $\approx1.8$ and the rate for the enhanced grid is $\approx3.3$. This is consistent with 
the theoretically results in Section~\ref{sec:SGError} which predicts the rate of the enhanced sparse grid to be approximately twice the rate
of the un-enhanced sparse grid.

Convergence of errors is shown with respect to increasing computational cost. Here, and throughout this paper, we define one unit of cost to be the
computational cost of one solve of the forward model and assume that the cost of solving the adjoint equation
is also equal to one unit.  For this example, we ignore the cost of computing the error estimate.  
We will challenge this last assumption shortly.

\subsubsection{Balancing physical and stochastic discretization errors}
\begin{figure}[ht]
\centering
\includegraphics[width=0.95\textwidth]{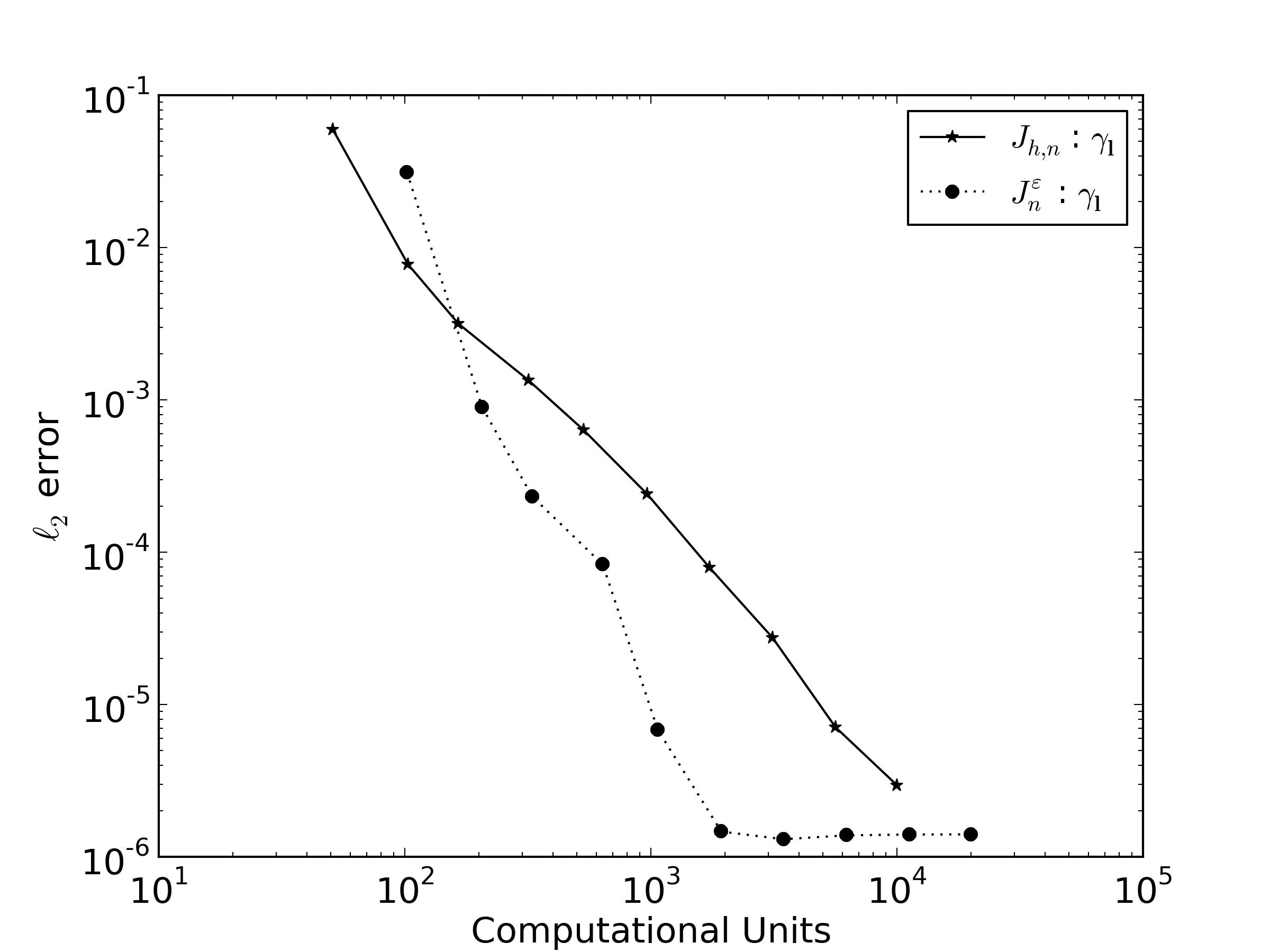}
\caption{Errors, in the Clenshaw-Curtis, dimension-adaptive, un-enhanced and enhanced sparse grid approximations of \eqref{eq:heat-eq-qoi} for the $25$d heat equation.}
\label{fig:heat-eq-convergence-d-25-unenhanced-vs-enahnced-brute-only}
\end{figure}

When quantifying the uncertainty of model predictions, both the deterministic and 
stochastic discretization errors must be accounted for.  It is inefficient to 
reduce the stochastic error to a level below the error introduced by the physical discretization.
This statement is supported by Figure~\ref{fig:heat-eq-convergence-d-25-unenhanced-vs-enahnced-brute-only}.
When when the stochastic error is smaller than the 
physical discretization error the accuracy in the enhanced sparse grid stagnates at approximately the 
physical discretization error. 

\subsubsection{Reducing the cost of computing a posterior error estimates}
Although typically ignored in the literature, the cost of computing $\eest$ is non-trivial. 
For this example the cost of evaluating $\eest$ is approximately $1/25$ of a computational unit
and the cost of evaluating $\eest$ at $100,000$ points, needed to compute the $\ell_2$ error, often outweighs the cost
of the forward and adjoint solves. 

\begin{figure}[ht]
\centering
\includegraphics[width=0.45\textwidth]{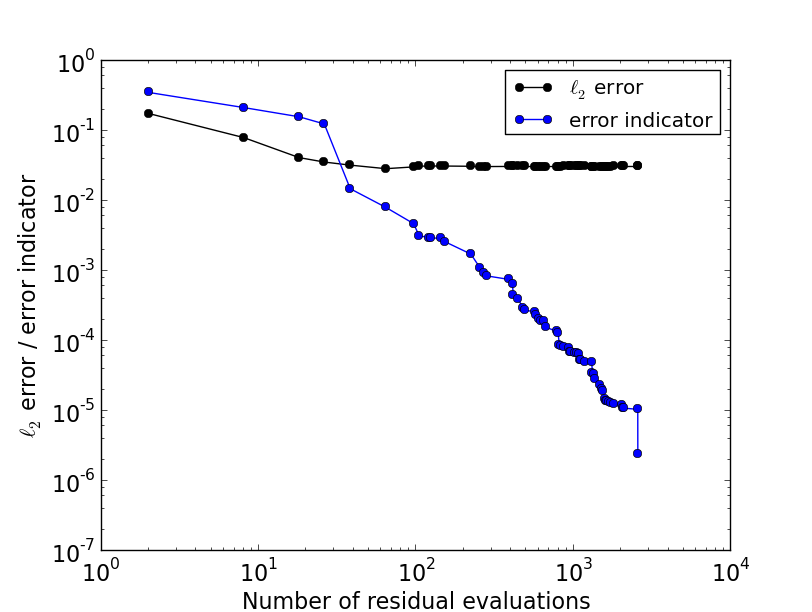}
\includegraphics[width=0.45\textwidth]{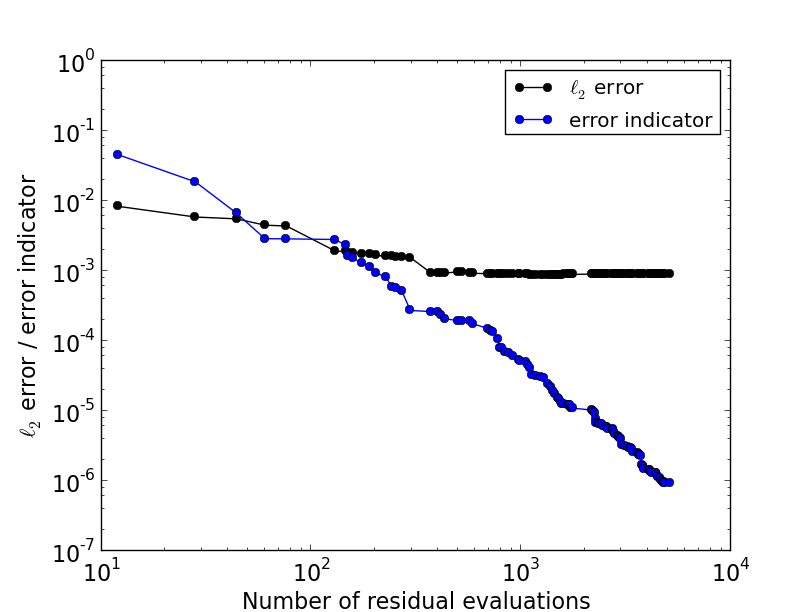}
\caption{The decay of the $\ell_2$ error vs the decay in the error indicator for the enhanced sparse grid $\qoienm$, starting from an initial grid
with 51 points (left) and 103 points (right).}
\label{fig:heat-eq-d-25-convergence-of-approximation-of-enhanced-sparse-grid}
\end{figure}

To reduce the cost of producing many a posterior error estimates we build an approximation $\qoienm$ of $\qoihn+\eest$ using the approach outlined in
Algorithm~\ref{alg:enhanced-interpolate}. Figure~\ref{fig:heat-eq-d-25-convergence-of-approximation-of-enhanced-sparse-grid} demonstrates that the
the performance of Algorithm~\ref{alg:enhanced-interpolate} is dependent on
the parameters $\tau$ and $n$. Simply driving the sparse grid error indicator
to zero is inefficient. After a certain number of samples, adding more points (each requiring a a posteriori error estimate) to the enhanced sparse grid has no effect 
on reducing the $\ell_2$ error of the approximation. However if $\tau$ and $n$ are set appropriately as described in~\ref{sec:error-balancing} the number of `wasted' samples can be minimized.

\begin{figure}[ht]
\centering
\includegraphics[width=0.95\textwidth]{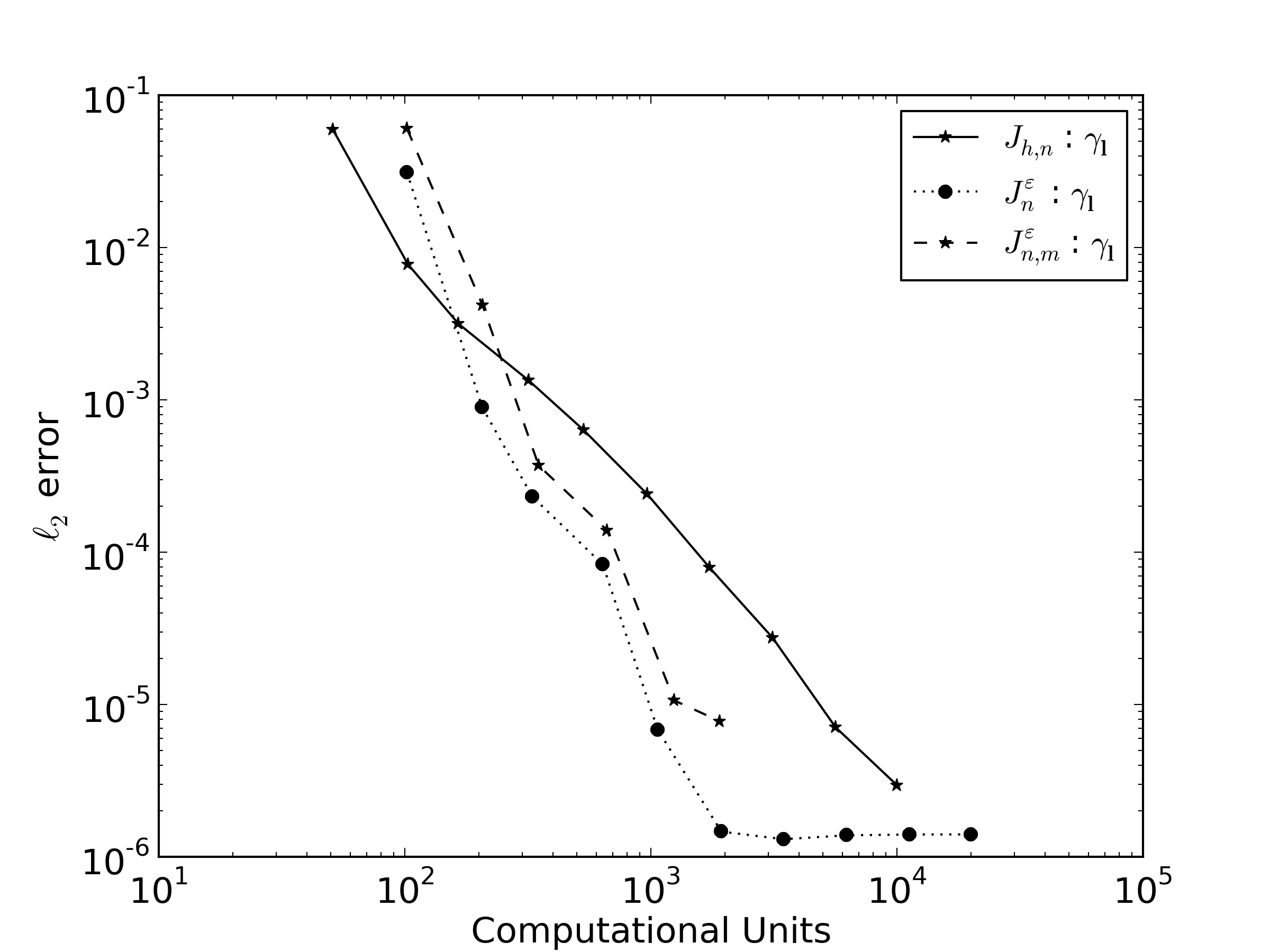}
\caption{Comparison of the $\ell_2$ convergence of the enhanced dimension-adaptive 
sparse grid using direct evaluation of the error estimate $\qoien$ 
and a sparse grid approximation of the enhanced QOI $\qoienm$.}
\label{fig:heat-eq-convergence-d-25-dim-adaptive}
\end{figure}

Figure~\ref{fig:heat-eq-convergence-d-25-dim-adaptive} plots the convergence of the 
un-enhanced surplus-based dimension-adaptive Clenshaw-Curtis sparse grids $\qoihn$,
$\qoien$ and $\qoienm$. The enhanced sparse grid $\qoienm$ is significantly more accurate than
the un-enhanced grid, however $\qoienm$ is not as accurate as $\qoien$. Building $\qoienm$ reduces the number of samples needed to generate
enhanced sparse grid samples at the expense of a slight degradation in accuracy. This degradation is acceptable as 
as in practice, adding $\eest$ to $\qoien$ is typically computationally infeasible\footnote{The computational cost of $\qoienm$ shown in Figure~\ref{fig:heat-eq-convergence-d-25-dim-adaptive} includes
the cost of evaluating the residual to compute the error estimates needed at the sparse
grid points, whereas we have assumed that computing $\eest$ for $\qoien$ has no cost.}.
From Figure~\ref{fig:heat-eq-convergence-d-25-dim-adaptive} we also observe that our choice of $\tau$ and $n$ stops refinement of $\qoienm$ before
the stochastic error is forced below the deterministic error $\delta_{\max}=8.41\times10^{-5}$. 

\begin{figure}[ht]
\centering
\includegraphics[width=0.95\textwidth]{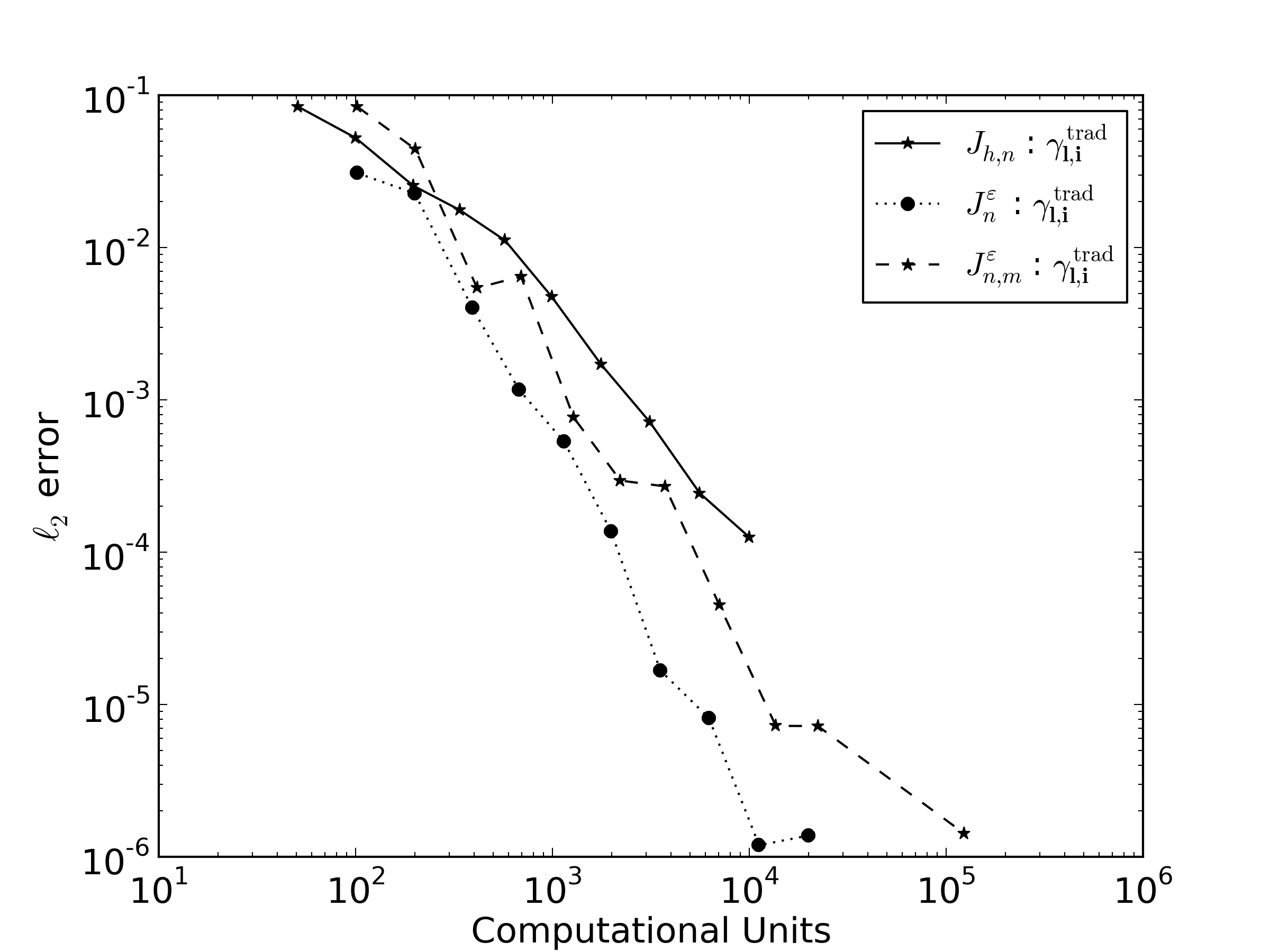}
\caption{Comparison of the $\ell_2$ convergence of the enhanced locally-adaptive sparse grid using direct evaluation 
of the error estimate and a sparse grid approximation of the enhanced QOI}
\label{fig:heat-eq-convergence-d-25-locally-adaptive}
\end{figure}

Figure~\ref{fig:heat-eq-convergence-d-25-locally-adaptive} plots the convergence of the 
un-enhanced surplus-based locally-adaptive sparse grids $\qoihn$,
$\qoien$ and $\qoienm$. The same conclusions drawn when using dimension adaptivity can 
also be made here. The improvement obtained by using $\qoienm$ instead of $\qoihn$ is reduced
however this can be addressed by the use of different refinement strategies discussed in 
Section~\ref{sec:results-refinement-strategies}

\subsection{Non-linear coupled system of ODEs}\label{sec:predator-prey}
For our second example, consider the non-linear system of ordinary differential equations governing a competitive Lotka–Volterra model of the 
population dynamics of species competing for some common resource. The model is given by 
\begin{equation}\label{eq:predprey}
\begin{cases}
\frac{d\sol_i}{dt} = r_i\sol_i\left(1-\sum_{j=1}^3\alpha_{ij}\sol_j\right), & \quad t\in (0,10],\\
\sol_i(0) = \sol_{i,0}
\end{cases},
\end{equation}
for $i = 1,2,3$.  The initial condition, $\sol_{i,0}$, and the self-interacting terms, $\alpha_{ii}$, are given, but the remaining interaction parameters, $\alpha_{ij}$ with $i\neq j$ as well as the re-productivity parameters, $r_i$, are unknown.
We assume that these parameters are each uniformly distributed on $[0.3,0.7]$.
We approximate the solution to~\eqref{eq:predprey} in time using a Backward Euler method with 1000 time steps ($\Delta t = 0.01$) which gives a deterministic error of 
approximately $1.00\times10^{-4}$.

The quantity of interest is the population of the third species at the final time, $\sol_3(10)$.  The corresponding adjoint problem is 
\begin{equation}\label{eq:adjpredprey}
\begin{cases}
-\frac{d\bphi_i}{dt} = r_i\bphi_i\left(1-\sum_{j=1}^3\alpha_{ij}\sol_j\right) + r_i\sol_i\left(1-\sum_{j=1}^3\alpha_{ji}\bphi_j\right), & \quad t\in (10,0],\\
\bphi_i(10) = 0, & i=1,2 \\
\bphi_3(10) = 1.
\end{cases},
\end{equation}
We approximate the adjoint solution in time using a second-order Crank-Nicholson method with the same number of time steps.

\begin{figure}[ht]
\centering
\includegraphics[width=0.95\textwidth]{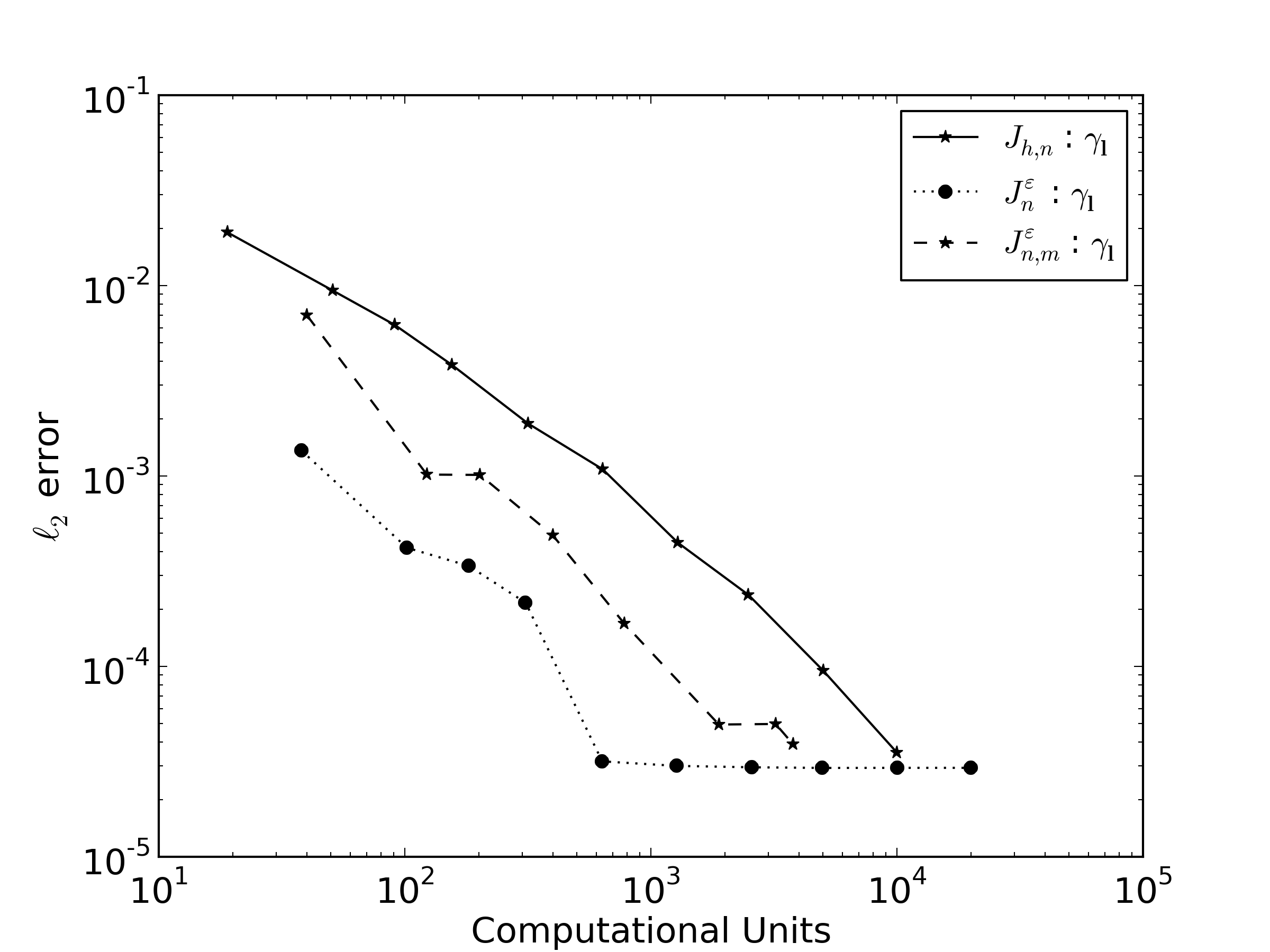}
\caption{Comparison of the $\ell_2$ convergence of the enhanced dimension-adaptive sparse grid using direct evaluation 
of the error estimate and a sparse grid approximation of the enhanced QOI}
\label{fig:predator-prey-convergence-d-9-dim-adaptive}
\end{figure}

\begin{figure}[ht]
\centering
\includegraphics[width=0.95\textwidth]{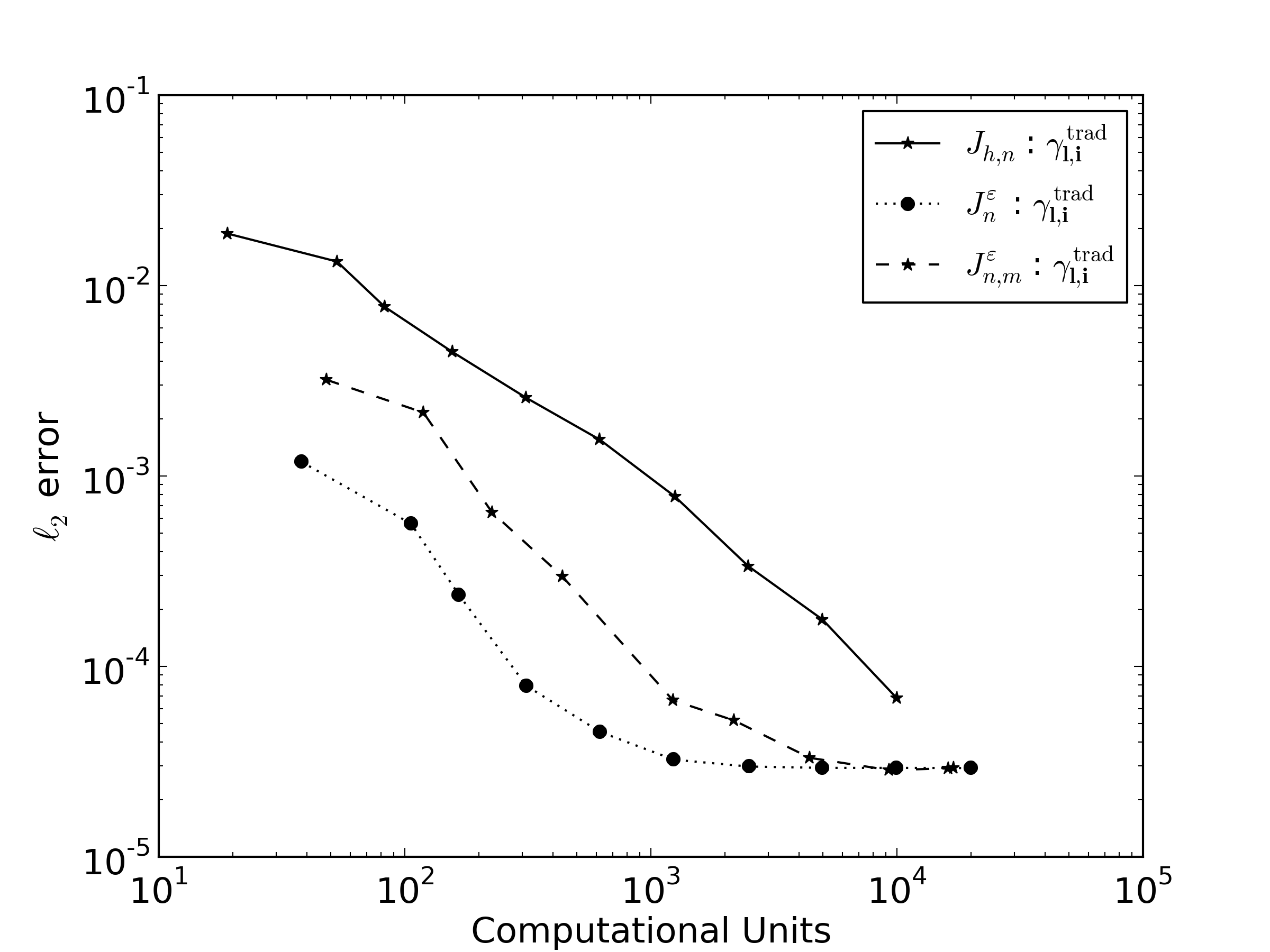}
\caption{Comparison of the $\ell_2$ convergence of the enhanced locally-adaptive sparse grid using direct evaluation 
of the error estimate and a sparse grid approximation of the enhanced QOI}
\label{fig:predator-prey-convergence-d-9-locally-adaptive}
\end{figure}

Figure~\ref{fig:predator-prey-convergence-d-9-dim-adaptive} and 
Figure~\ref{fig:predator-prey-convergence-d-9-locally-adaptive} 
compare the convergence
of $\qoihn$, $\qoien$ and $\qoienm$ when using a surplus based 
dimension-adaptive Clenshaw-Curtis sparse grid and surplus locally adaptive sparse
grid, respectively. Again, significant increases in accuracy 
can be achieved by using $\qoienm$ instead of $\qoihn$.  In addition we again see that we are able to correctly stop refinement approximately 
when the stochastic error of the sparse grid approximation is the same order as the physical discretization error.
We emphasize that, as in the previous example, 
the accuracy of $\qoienm$ can be forced towards the accuracy of the best possible enhanced approximation $\qoien+\eest$ at the cost of
additional residual calculations (error estimates) at new sparse grid points.

\subsection{A posteriori error based refinement strategies}\label{sec:results-refinement-strategies}
In this section we demonstrate the utility of using a posteriori error based refinement strategies
outlined in Section~\ref{sec:sparse-grid-adaptivity}. Specifically we consider a posteriori 
refinement indicators for both dimension and locally adapted sparse grid. In the 
case of local refinement we also consider the simultaneous dimension and local refinement
in addition to the typical local refinement strategy used in all examples thus far. 
We remark that in all the figures shown in the remainder of this paper, the cost of computing the sparse grid approximation, shown on the horizontal axis,
includes the cost of evaluating the a posteriori error estimate, building the enhanced the sparse grid approximation, and calculating the refinement
indicators.
\subsubsection{Dimension refinement}
\begin{figure}[ht]
\centering
\includegraphics[width=0.95\textwidth]{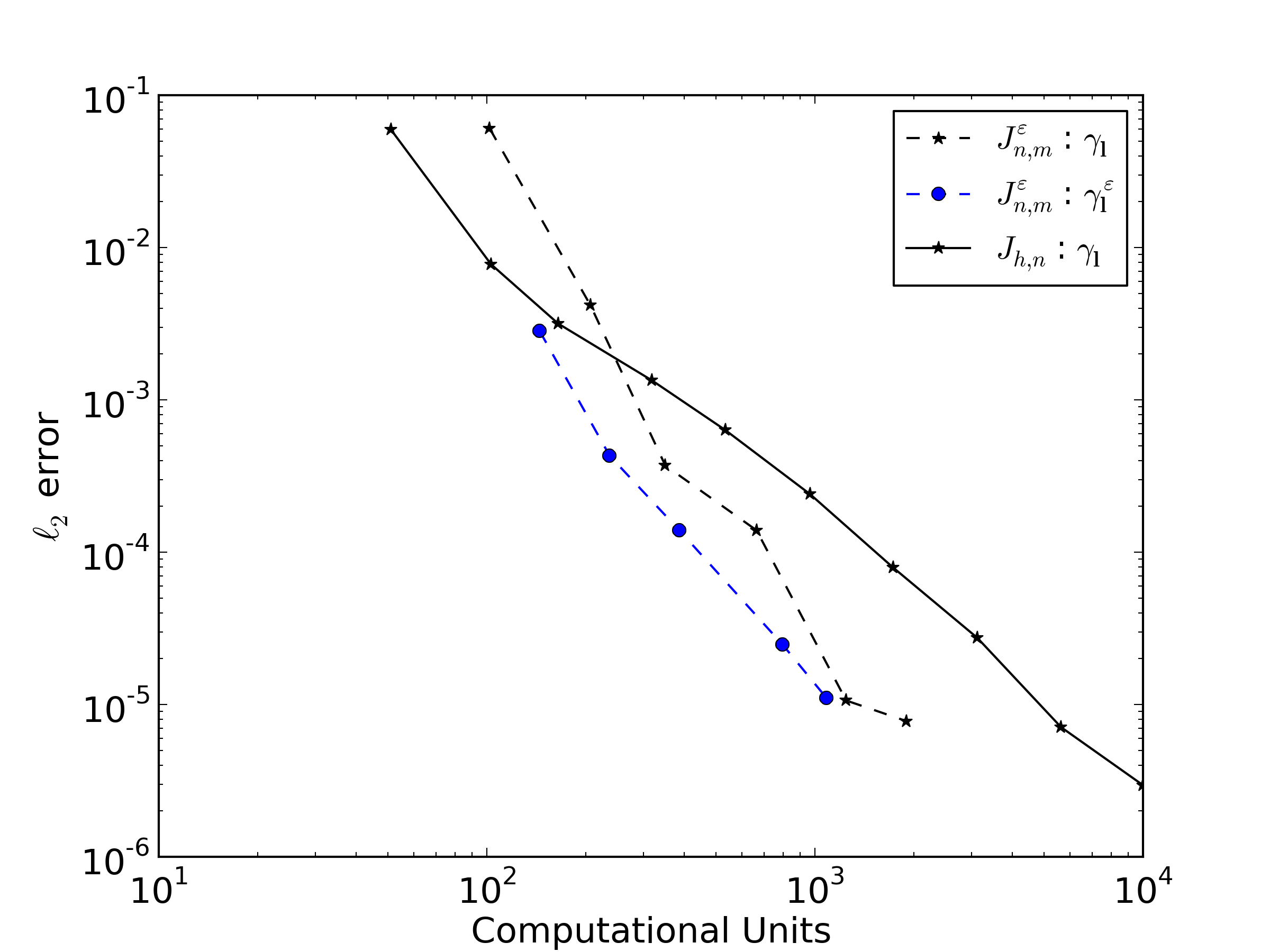}
\caption{Comparison of dimension-adaptive refinement strategies when applied to the 25d
heat equation}
\label{fig:25-d-heat-equation-dim-adaptive-refinement-criteria-comparision}
\end{figure}

Figure~\ref{fig:25-d-heat-equation-dim-adaptive-refinement-criteria-comparision}
compares refinement strategies for the 25d heat equation problem presented in 
Section~\ref{sec:heat-equation}. It is clear that the use of a posteriori refinement
results in significant increases in efficiency. This result is due to the fact that
the a posteriori subspace refinement indicator $\gamma_\bl$ only requires computation of error
estimates which are relatively cheap compared to the evaluation of the forward model. 
The function evaluations that are needed by surplus refinement to probe candidate 
subspaces are often redundant and this redundant evaluation can be avoided by using
a posteriori refinement. 

\begin{figure}[ht]
\centering
\includegraphics[width=0.95\textwidth]{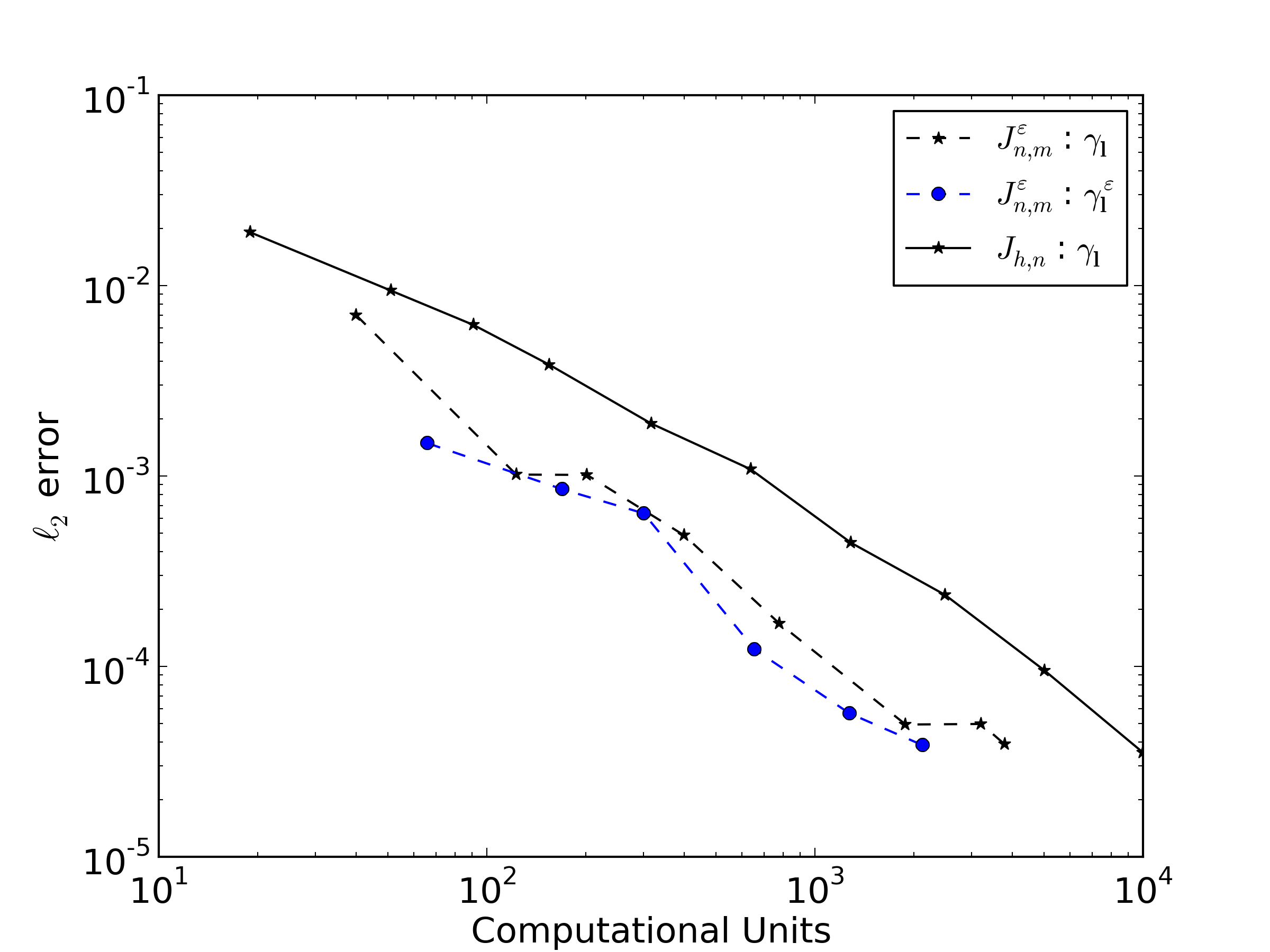}
\caption{Comparison of dimension-adaptive refinement strategies when applied to 
9d predator prey model}
\label{fig:9d-predator-prey-dim-adaptive-refinement-criteria-comparision}
\end{figure}

The high-dimensional nature and the strong degree of anisotropy of this 
problem are an ideal setting for a posteriori refinement. If the function was lower dimensional
or have less anisotropy then surplus refinement would be more competitive as the amount
of redundant model evaluations would decrease. This point is supported by
Figure~\ref{fig:9d-predator-prey-dim-adaptive-refinement-criteria-comparision} which
compares a posteriori and surplus dimension adaptive sparse grids when used to 
solve the 9d predator-prey model presented in Section~\ref{sec:predator-prey}. The predator prey model
is less anisotropic than the diffusion model and thus the benefits of a-posteriori refinement are reduced, although the
benefits are still positive.

\subsubsection{Local refinement}
\begin{figure}[ht]
\centering
\includegraphics[width=0.95\textwidth]{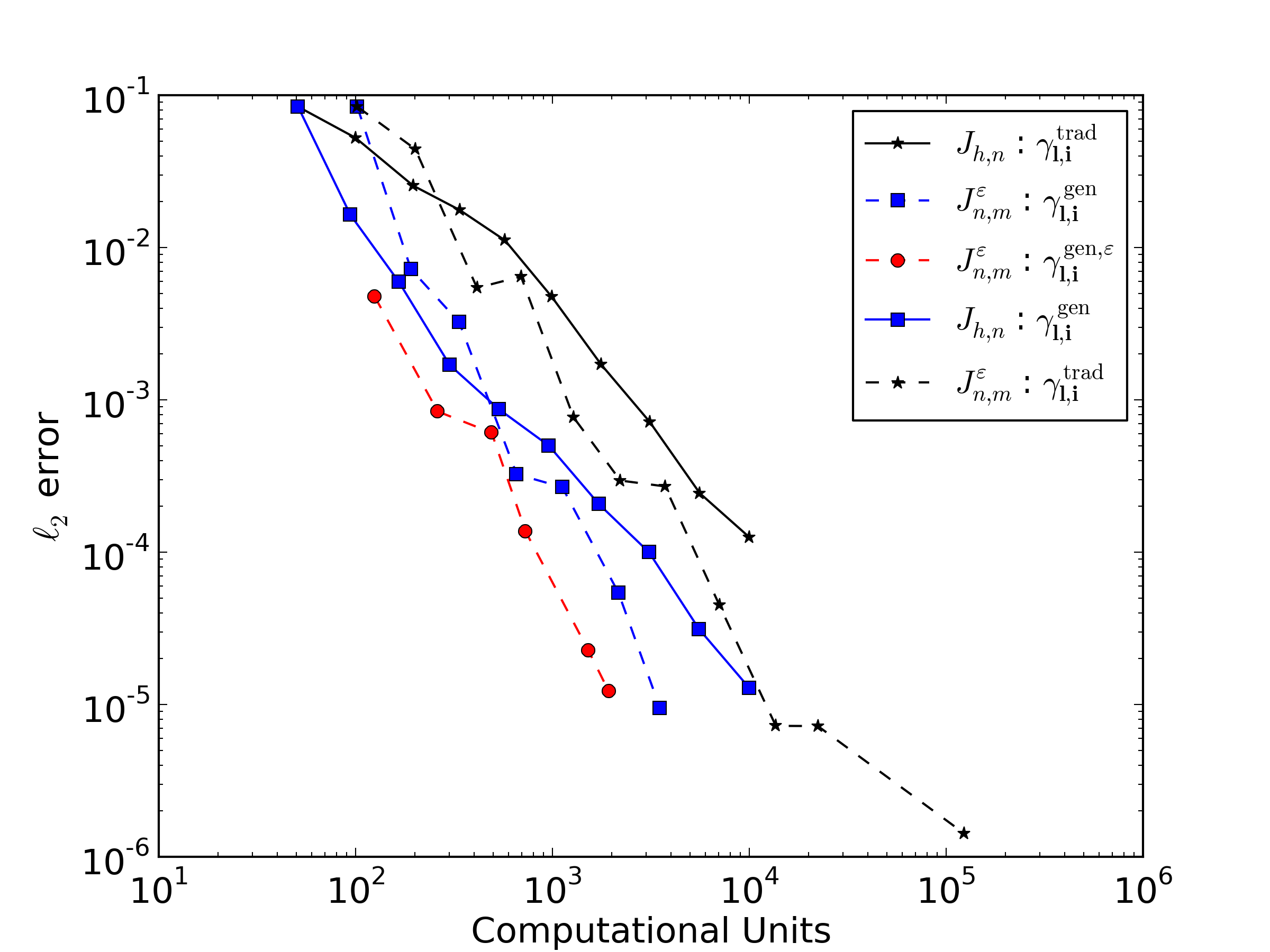}
\caption{Comparison of local refinement strategies when applied to the 25d heat equation}
\label{fig:25d-heat-equation-locally-adaptive-refinement-criteria-comparision}
\end{figure}

A posteriori refinement can also be used in conjunction with local refinement. 
Figure~\ref{fig:25d-heat-equation-locally-adaptive-refinement-criteria-comparision}
compares refinement strategies for the 25d heat equation. The use of generalized local
refinement results in a vast improvement over traditional local refinement. The best refinement
strategy for this problem is to use a posteriori refinement indicators with the generalized local
refinement. This result is due to the high-dimension of the problem and the high degree 
of anisotropy.

\begin{figure}[ht]
\centering
\includegraphics[width=0.95\textwidth]{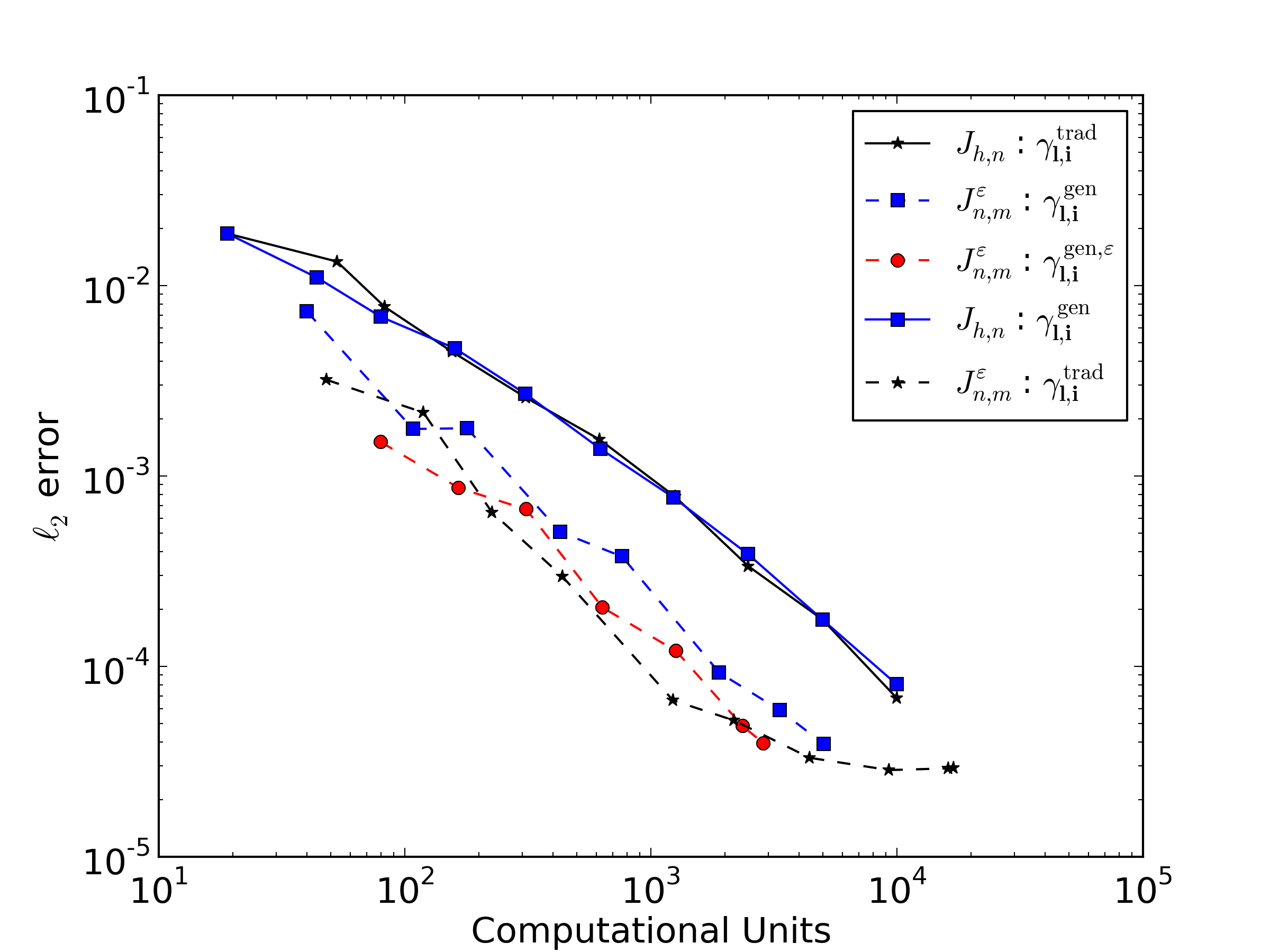}
\caption{Comparison of local refinement strategies when applied to the 
9d predator prey model}
\label{fig:9d-predator-prey-locally-adaptive-refinement-criteria-comparision}
\end{figure}

Figure~\ref{fig:9d-predator-prey-locally-adaptive-refinement-criteria-comparision}
compares refinement strategies for the 9d predator prey model. Due to the low degree
of anisotropy of this function generalized local refinement provides no improvement 
over traditional local refinement. In the case of the surplus-based enhanced sparse grid 
the accuracy is slightly worse. Traditional surplus based refinement is the most accurate
strategy here, but the a posteriori error based generalized local refinement has comparable
accuracy.

The best refinement strategy is problem dependent. The authors recommend basing the choice of refinement strategy by performing studies such as the one presented in~\cite{winkor2013}. However
in the absence of such a investigation and based upon the examples presented in this paper, the authors recommend using a posteriori based refinement for both dimension adaptive and local adaptive sparse
grids. When using local adaptivity a posteriori based refinement should be combined with the generalized local adaptive procedure presented in~\cite{jakeman2013localuq}.

Finally we remark that the primary concern of this section was to discuss the efficacy of a posteriori and surplus based refinement and not the strengths and weaknesses
of dimension vs local refinement. Examples using both dimension and local refinement were provided solely to demonstrate that either type of refinement can be used
in conjunction with a posteriori error estimates. The benefit of local and dimension refinement is dependent on the regularity of the PDE solutions. 
Practically determining which approach is better is an open question.
\section{Conclusions}\label{sec:conclusions}
In this paper we present an algorithm for adaptive sparse grid approximations
of quantities of interest computed from discretized partial differential equations.
We use adjoint-based a posteriori error estimates of the 
interpolation error in the sparse grid to enhance the sparse grid approximation. 
Using a number of numerical examples we show that utilizing a posteriori error estimates provides significantly 
more accurate functional values for random samples of the sparse grid approximation. The cost of computing these
error estimates can be non-trivial and thus we provide a practical method for enhancing a sparse grid approximation
using only a finite set of error estimates. 

Aside from using a posteriori error estimates to enhance an approximation we also demonstrate that error estimates can be used to increase
the efficiency of adaptive sampling. Traditional sparse grid adaptivity employs error indicators based upon the hierarchical surplus 
are used to flag dimensions or local regions for refinement.  However such 
approaches require the model to be evaluated at a new point before one can 
determine if refinement should have taken place. In this paper we numerically demonstrated that refinement using a posteriori error estimates
can significantly reduce the amount of redundant sampling compared when compared to traditional hierarchical refinement. Refinement indicators based 
upon a posteriori error estimates are most effective when the quantity of interest being approximated is anisotropic, however even in the absence
of anisotropy a posteriori indicators perform at least comparably to traditional hierarchical surplus based indicators.

In combination with the aforementioned enhancement and refinement procedures we use a posteriori error estimates 
to ensure that the sparse grid is not refined beyond the point at which the stochastic interpolation error is below the physical discretization
error. The methodology presented provides a practical means of balancing the stochastic and deterministic discretization errors.

\bibliographystyle{acm}
\bibliography{error-refs}

\end{document}